\pdfoutput=1
\documentclass[%
  sigconf,%
  twocolumn,pbalance=true,
  nonacm=true,
]{acmart}

\usepackage[utf8]{inputenc}
\usepackage{amsmath}
\usepackage{mathtools}
\usepackage{algorithm}
\usepackage{algorithmicx}
\usepackage{algpseudocode}

\usepackage{faktor}

\usepackage{enumerate}

\usepackage{pgf,pgfcore,pgfkeys}
\usepackage{tikz}
\usepackage{tikz-cd}
\usepackage{tikz,tkz-euclide,graphicx}
\usetikzlibrary{through}
\usetikzlibrary{arrows}
\usetikzlibrary{fit,matrix,positioning,svg.path}

\usepackage[]{xcolor}
\definecolor{orange}{rgb}{0.85,0.4,.1}
\definecolor{burgundy}{rgb}{0.65,0.0,0}
\definecolor{blueb}{rgb}{0.0,0.3,0.5}

\usepackage[capitalize,noabbrev]{cleveref}
\usepackage{setspace}

\usepackage{subfigure}

\makeatletter
\newcounter{algorithmicH}
\let\oldalgorithmic\algorithmic
\renewcommand{\algorithmic}{%
  \stepcounter{algorithmicH}
  \oldalgorithmic}
\renewcommand{\theHALG@line}{ALG@line.\thealgorithmicH.\arabic{ALG@line}}
\makeatother

\newcommand{\N}{\mathbb{N}}
\newcommand{\Z}{\mathbb{Z}}
\newcommand{\Q}{\mathbb{Q}}

\newcommand{\CC}{\mathbb{C}}
\newcommand{\F}{\mathbb{F}}

\newcommand{\Norm}[1]{\left\Vert #1 \right\Vert}
\newcommand{\Ninf}[1]{\Norm{#1}_{\infty}}
\newcommand{\NI}[1]{\Norm{x}_{1}}

\newcommand{\abs}[1]{\lvert #1 \rvert}

\newcommand{\norm}[1]{\mathrm{N}_{#1}}

\newcommand{\Gal}[2]{\mathrm{Gal}(#1/#2)}

\newcommand{\CRT}[1]{$\mathrm{CRT}$}

\newcommand{\U}[1]{\mathcal{O}_{#1}^\times}

\renewcommand{\O}[1]{\mathcal{O}_{#1}}

\newcommand{\segment}[2]{\llbracket #1, #2 \rrbracket}

\newcommand{\mf}[1]{\mathfrak{#1}}

\newcommand{\rd}[2]{f(#1 \vert #2)}

\newcommand{\Cinf}{\mathsf{C}_{\infty}}
\newcommand{\Mpol}[2]{\mathsf{M}(#1,#2)}
\newcommand{\Mint}[1]{\mathsf{M}(#1)}
\newcommand{\FFroot}[2]{\mathsf{R}(#1,#2)}
\newcommand{\llog}{\log\!\log}

\DeclareMathOperator{\lcm}{lcm}

\theoremstyle{plain}
\newtheorem{theorem}{Theorem}
\newtheorem{proposition}{Proposition}
\newtheorem{lemma}{Lemma}

\theoremstyle{definition}

\newtheorem{remark}{Remark}

\newtheorem*{notation}{Notation}

\newcommand{\ie}{i.e.,~}
\newcommand{\eg}{e.g.,~}
\newcommand{\wrt}{w.r.t.~}
\newcommand{\resp}{resp.~}

\usepackage{scalerel}
\definecolor{orcidlogocol}{HTML}{A6CE39}
\tikzset{
    orcidlogo/.pic={
        \fill[orcidlogocol] svg{M256,128c0,70.7-57.3,128-128,128C57.3,256,0,198.7,0,128C0,57.3,57.3,0,128,0C198.7,0,256,57.3,256,128z};
        \fill[white] svg{M86.3,186.2H70.9V79.1h15.4v48.4V186.2z}
                     svg{M108.9,79.1h41.6c39.6,0,57,28.3,57,53.6c0,27.5-21.5,53.6-56.8,53.6h-41.8V79.1z M124.3,172.4h24.5c34.9,0,42.9-26.5,42.9-39.7c0-21.5-13.7-39.7-43.7-39.7h-23.7V172.4z}
                     svg{M88.7,56.8c0,5.5-4.5,10.1-10.1,10.1c-5.6,0-10.1-4.6-10.1-10.1c0-5.6,4.5-10.1,10.1-10.1C84.2,46.7,88.7,51.3,88.7,56.8z};
    }
}

\makeatletter
\DeclareRobustCommand{\orcid@link}[1]{
    \href{https://orcid.org/#1}{%
        \mbox{\scalerel*{
            \begin{tikzpicture}[yscale=-1,transform shape]
                \pic{orcidlogo};
            \end{tikzpicture}
        }{|}}%
    }%
}
\newcommand{\orcidlink}[1]{\texorpdfstring{\unexpanded{\orcid@link{#1}}}{}}
\makeatother

\newcommand{\giturl}{https://github.com/ob3rnard/eth-roots}
\newcommand{\olivierOrcID}{0000-0002-3410-3134}
\newcommand{\pierrealainOrcID}{0000-0003-4997-2276}
\newcommand{\andreaOrcID}{0000-0001-8318-4922}

\title{Computing \texorpdfstring{\(e\)}{e}-th roots in number fields}

\author{Olivier Bernard\orcidlink{\olivierOrcID}}
\email{olivier.bernard@normalesup.org}
\orcid{\olivierOrcID}
\affiliation{%
  \institution{Zama}%
  \city{Paris}%
  \country{France}%
}

\author{Pierre-Alain Fouque\orcidlink{\pierrealainOrcID}}
\email{pierre-alain.fouque@irisa.fr}
\orcid{\pierrealainOrcID}
\affiliation{%
  \institution{Univ Rennes, IRISA}%
  \city{Rennes}%
  \country{France}%
}

\author{Andrea Lesavourey\orcidlink{\andreaOrcID}}
\email{andrea.lesavourey@irisa.fr}
\orcid{\andreaOrcID}
\affiliation{%
  \institution{Univ Rennes, CNRS, IRISA}
  \city{Rennes}
  \country{France}
}

\newcommand{\acknowledgements}{%
  The authors would like to thank Dr.~Razvan Barbulescu and Dr.~Aurel Page for suggesting using Schirokauer maps to detect huge $e$-th powers.
  The authors would also like to thank the anonymous reviewers for their interesting and useful comments.
  Most of this work was done while Olivier Bernard was employed by Thales, Gennevilliers, France.
  Andrea Lesavourey is funded by the Direction Générale de l’Armement (Pôle de Recherche CYBER), with the support of Région Bretagne.%
}
\thanks{\acknowledgements}

\copyrightyear{2023}
\IfFileExists{root_crt.aux}{}{%
  \IfFileExists{root_crt_aux.tex}{\makeatletter\relax 
\providecommand\hyper@newdestlabel[2]{}
\providecommand\zref@newlabel[2]{}
\providecommand\HyperFirstAtBeginDocument{\AtBeginDocument}
\HyperFirstAtBeginDocument{\ifx\hyper@anchor\@undefined
\global\let\oldcontentsline\contentsline
\gdef\contentsline#1#2#3#4{\oldcontentsline{#1}{#2}{#3}}
\global\let\oldnewlabel\newlabel
\gdef\newlabel#1#2{\newlabelxx{#1}#2}
\gdef\newlabelxx#1#2#3#4#5#6{\oldnewlabel{#1}{{#2}{#3}}}
\AtEndDocument{\ifx\hyper@anchor\@undefined
\let\contentsline\oldcontentsline
\let\newlabel\oldnewlabel
\fi}
\fi}
\global\let\hyper@last\relax 
\gdef\HyperFirstAtBeginDocument#1{#1}
\providecommand\HyField@AuxAddToFields[1]{}
\providecommand\HyField@AuxAddToCoFields[2]{}
\citation{number_field_sieve,couveignes}
\citation{norm_relations}
\citation{tw-sti}
\citation{tw-sti}
\citation{thome}
\citation{belabas_roots}
\citation{couveignes}
\citation{sagemath}
\@writefile{toc}{\contentsline {section}{Abstract}{1}{section*.1}\protected@file@percent }
\@writefile{toc}{\contentsline {section}{\numberline {1}Introduction}{1}{section.1}\protected@file@percent }
\newlabel{sec:introduction}{{1}{1}{Introduction}{section.1}{}}
\newlabel{sec:introduction@cref}{{[section][1][]1}{[1][1][]1}}
\citation{PARI2}
\citation{tw-sti}
\citation{GG13}
\citation{belabas_roots}
\citation{fieker-reconstruction}
\citation{AMM77}
\citation{cohen-course}
\citation{Sho93}
\citation{thome}
\citation{belabas_roots}
\@writefile{toc}{\contentsline {section}{\numberline {2}Preliminaries}{2}{section.2}\protected@file@percent }
\newlabel{sec:preliminaries}{{2}{2}{Preliminaries}{section.2}{}}
\newlabel{sec:preliminaries@cref}{{[section][2][]2}{[1][2][]2}}
\@writefile{toc}{\contentsline {paragraph}{Complexities.}{2}{section*.6}\protected@file@percent }
\@writefile{toc}{\contentsline {subsection}{\numberline {2.1}Bounds on root coefficients}{2}{subsection.2.1}\protected@file@percent }
\newlabel{sec:bound}{{2.1}{2}{Bounds on root coefficients}{subsection.2.1}{}}
\newlabel{sec:bound@cref}{{[subsection][1][2]2.1}{[1][2][]2}}
\newlabel{lem:boundcoef}{{1}{2}{}{lemma.1}{}}
\newlabel{lem:boundcoef@cref}{{[lemma][1][]1}{[1][2][]2}}
\newlabel{lem:boundethroot}{{2}{2}{}{lemma.2}{}}
\newlabel{lem:boundethroot@cref}{{[lemma][2][]2}{[1][2][]2}}
\@writefile{toc}{\contentsline {subsection}{\numberline {2.2}Computing $e$-th roots in finite fields}{2}{subsection.2.2}\protected@file@percent }
\newlabel{ethfp}{{2.2}{2}{Computing \texorpdfstring {$e$}{e}-th roots in finite fields}{subsection.2.2}{}}
\newlabel{ethfp@cref}{{[subsection][2][2]2.2}{[1][2][]2}}
\@writefile{toc}{\contentsline {section}{\numberline {3}Generic local-global methods}{2}{section.3}\protected@file@percent }
\newlabel{sec:generic}{{3}{2}{Generic local-global methods}{section.3}{}}
\newlabel{sec:generic@cref}{{[section][3][]3}{[1][2][]2}}
\@writefile{toc}{\contentsline {subsection}{\numberline {3.1}$p$-adic lifting}{2}{subsection.3.1}\protected@file@percent }
\newlabel{padiclifting}{{3.1}{2}{\texorpdfstring {$p$}{p}-adic lifting}{subsection.3.1}{}}
\newlabel{padiclifting@cref}{{[subsection][1][3]3.1}{[1][2][]2}}
\citation{Bre75}
\citation{belabas_roots}
\citation{belabas_roots}
\citation{belabas_roots}
\citation{belabas_roots}
\citation{belabas_roots}
\citation{belabas_roots}
\citation{GN08}
\citation{thome}
\newlabel{newtiteration}{{1}{3}{\texorpdfstring {$p$}{p}-adic lifting}{equation.3.1}{}}
\newlabel{newtiteration@cref}{{[equation][1][]1}{[1][3][]3}}
\@writefile{loa}{\contentsline {algorithm}{\numberline {1}{\ignorespaces $p$-adic lifting for $e$-th root\relax }}{3}{algorithm.1}\protected@file@percent }
\providecommand*\caption@xref[2]{\@setref\relax\@undefined{#1}}
\newlabel{alg:padic}{{1}{3}{$p$-adic lifting for $e$-th root\relax }{algorithm.1}{}}
\newlabel{alg:padic@cref}{{[algorithm][1][]1}{[1][3][]3}}
\@writefile{toc}{\contentsline {subsection}{\numberline {3.2}$\mathfrak  {p}$-adic reconstruction}{3}{subsection.3.2}\protected@file@percent }
\newlabel{padicreconstruction}{{3.2}{3}{\texorpdfstring {$\mf {p}$}{P}-adic reconstruction}{subsection.3.2}{}}
\newlabel{padicreconstruction@cref}{{[subsection][2][3]3.2}{[1][3][]3}}
\@writefile{toc}{\contentsline {section}{\numberline {4}Using Chinese Remainder Theorem: the easy cases}{3}{section.4}\protected@file@percent }
\newlabel{sec:roots-easy}{{4}{3}{Using Chinese Remainder Theorem: the easy cases}{section.4}{}}
\newlabel{sec:roots-easy@cref}{{[section][4][]4}{[1][3][]3}}
\newlabel{eq:goodcase}{{2}{3}{Using Chinese Remainder Theorem: the easy cases}{equation.4.2}{}}
\newlabel{eq:goodcase@cref}{{[equation][2][]2}{[1][3][]3}}
\@writefile{toc}{\contentsline {subsection}{\numberline {4.1}A CRT-based method for $e$-th roots}{3}{subsection.4.1}\protected@file@percent }
\citation{GG13}
\citation{GG13}
\citation{GG13}
\citation{GG13}
\citation{Molin}
\citation{neukirch}
\citation{neukirch}
\citation{tw-sti}
\citation{tw-sti}
\citation{trager}
\citation{fieker-reconstruction,belabas_roots}
\@writefile{loa}{\contentsline {algorithm}{\numberline {2}{\ignorespaces Number field CRT for $e$-th root mod $q$\relax }}{4}{algorithm.2}\protected@file@percent }
\newlabel{alg:crt-primes}{{2}{4}{Number field CRT for $e$-th root mod $q$\relax }{algorithm.2}{}}
\newlabel{alg:crt-primes@cref}{{[algorithm][2][]2}{[1][3][]4}}
\newlabel{nfcrt:stepmodid}{{5}{4}{Number field CRT for $e$-th root mod $q$\relax }{ALG@line.ALG@line.2.5}{}}
\newlabel{nfcrt:stepmodid@cref}{{[line][5][]5}{[1][3][]4}}
\newlabel{nfcrt:idcrt}{{7}{4}{Number field CRT for $e$-th root mod $q$\relax }{ALG@line.ALG@line.2.7}{}}
\newlabel{nfcrt:idcrt@cref}{{[line][7][]7}{[1][3][]4}}
\@writefile{loa}{\contentsline {algorithm}{\numberline {3}{\ignorespaces Double CRT for $e$-th root\relax }}{4}{algorithm.3}\protected@file@percent }
\newlabel{alg:crt-root}{{3}{4}{Double CRT for $e$-th root\relax }{algorithm.3}{}}
\newlabel{alg:crt-root@cref}{{[algorithm][3][]3}{[1][3][]4}}
\newlabel{dblcrt:stepmodq}{{5}{4}{Double CRT for $e$-th root\relax }{ALG@line.ALG@line.3.5}{}}
\newlabel{dblcrt:stepmodq@cref}{{[line][5][]5}{[1][3][]4}}
\newlabel{dblcrt:lastcrt}{{7}{4}{Double CRT for $e$-th root\relax }{ALG@line.ALG@line.3.7}{}}
\newlabel{dblcrt:lastcrt@cref}{{[line][7][]7}{[1][3][]4}}
\newlabel{pr:cplxcrtq}{{2}{4}{}{proposition.2}{}}
\newlabel{pr:cplxcrtq@cref}{{[proposition][2][]2}{[1][3][]4}}
\newlabel{rem:crt-primes}{{1}{4}{}{remark.1}{}}
\newlabel{rem:crt-primes@cref}{{[remark][1][]1}{[1][4][]4}}
\@writefile{toc}{\contentsline {subsection}{\numberline {4.2}Bad cases: existence of cyclotomic subfield}{4}{subsection.4.2}\protected@file@percent }
\newlabel{thm:equiv-bad-cases}{{1}{4}{}{theorem.1}{}}
\newlabel{thm:equiv-bad-cases@cref}{{[theorem][1][]1}{[1][4][]4}}
\newlabel{lem:neuk}{{3}{4}{Bauer in~\cite {neukirch}}{lemma.3}{}}
\newlabel{lem:neuk@cref}{{[lemma][3][]3}{[1][4][]4}}
\@writefile{toc}{\contentsline {subsection}{\numberline {4.3}Experimental results}{4}{subsection.4.3}\protected@file@percent }
\newlabel{sec:thome-experimental-results}{{4.3}{4}{Experimental results}{subsection.4.3}{}}
\newlabel{sec:thome-experimental-results@cref}{{[subsection][3][4]4.3}{[1][4][]4}}
\citation{tw-sti}
\citation{couveignes,thome}
\newlabel{subfig:crt-3}{{1(a)}{5}{Subfigure 1(a)}{subfigure.1.1}{}}
\newlabel{subfig:crt-3@cref}{{[subfigure][1][1]1(a)}{[1][4][]5}}
\newlabel{sub@subfig:crt-3}{{(a)}{5}{Subfigure 1(a)\relax }{subfigure.1.1}{}}
\newlabel{subfig:crt-71}{{1(b)}{5}{Subfigure 1(b)}{subfigure.1.2}{}}
\newlabel{subfig:crt-71@cref}{{[subfigure][2][1]1(b)}{[1][4][]5}}
\newlabel{sub@subfig:crt-71}{{(b)}{5}{Subfigure 1(b)\relax }{subfigure.1.2}{}}
\newlabel{subfig:crt-1637}{{1(c)}{5}{Subfigure 1(c)}{subfigure.1.3}{}}
\newlabel{subfig:crt-1637@cref}{{[subfigure][3][1]1(c)}{[1][4][]5}}
\newlabel{sub@subfig:crt-1637}{{(c)}{5}{Subfigure 1(c)\relax }{subfigure.1.3}{}}
\newlabel{subfig:crt-13099}{{1(d)}{5}{Subfigure 1(d)}{subfigure.1.4}{}}
\newlabel{subfig:crt-13099@cref}{{[subfigure][4][1]1(d)}{[1][4][]5}}
\newlabel{sub@subfig:crt-13099}{{(d)}{5}{Subfigure 1(d)\relax }{subfigure.1.4}{}}
\@writefile{lof}{\contentsline {figure}{\numberline {1}{\ignorespaces Timings (s) for \texttt  {nfroots} and \Cref  {alg:crt-root} plotted against the degree, for various prime $e$ over cyclotomic fields.\relax }}{5}{figure.caption.7}\protected@file@percent }
\newlabel{fig:thome-naive-root}{{1}{5}{Timings (s) for \texttt {nfroots} and \Cref {alg:crt-root} plotted against the degree, for various prime $e$ over cyclotomic fields.\relax }{figure.caption.7}{}}
\newlabel{fig:thome-naive-root@cref}{{[figure][1][]1}{[1][4][]5}}
\@writefile{lof}{\contentsline {subfigure}{\numberline{(a)}{\ignorespaces {\(e = 3\)}}}{5}{figure.caption.7}\protected@file@percent }
\@writefile{lof}{\contentsline {subfigure}{\numberline{(b)}{\ignorespaces {\(e = 71\)}}}{5}{figure.caption.7}\protected@file@percent }
\@writefile{lof}{\contentsline {subfigure}{\numberline{(c)}{\ignorespaces {\(e = 1637\)}}}{5}{figure.caption.7}\protected@file@percent }
\@writefile{lof}{\contentsline {subfigure}{\numberline{(d)}{\ignorespaces {\(e = 13099\)}}}{5}{figure.caption.7}\protected@file@percent }
\@writefile{lot}{\contentsline {table}{\numberline {1}{\ignorespaces Timings (s) for \(e\)-th roots within the saturation of Stickelberger \(S\)-units \citep  {tw-sti} for selected cyclotomic fields\relax }}{5}{table.caption.8}\protected@file@percent }
\newlabel{tab:timings-s-units}{{1}{5}{Timings (s) for \(e\)-th roots within the saturation of Stickelberger \(S\)-units \cite {tw-sti} for selected cyclotomic fields\relax }{table.caption.8}{}}
\newlabel{tab:timings-s-units@cref}{{[table][1][]1}{[1][4][]5}}
\@writefile{toc}{\contentsline {section}{\numberline {5}A relative Couveignes' method: the bad cases}{5}{section.5}\protected@file@percent }
\newlabel{sec:roots-bad}{{5}{5}{A relative Couveignes' method: the bad cases}{section.5}{}}
\newlabel{sec:roots-bad@cref}{{[section][5][]5}{[1][5][]5}}
\@writefile{toc}{\contentsline {subsection}{\numberline {5.1}A relative Couveignes' method}{5}{subsection.5.1}\protected@file@percent }
\newlabel{sec:relat-couv}{{5.1}{5}{A relative Couveignes' method}{subsection.5.1}{}}
\newlabel{sec:relat-couv@cref}{{[subsection][1][5]5.1}{[1][5][]5}}
\newlabel{lem:norm-com-inert}{{4}{5}{}{lemma.4}{}}
\newlabel{lem:norm-com-inert@cref}{{[lemma][4][]4}{[1][5][]5}}
\newlabel{lem:norm-bij}{{5}{5}{}{lemma.5}{}}
\newlabel{lem:norm-bij@cref}{{[lemma][5][]5}{[1][5][]5}}
\citation{cohen-advanced,neukirch}
\citation{lang2012algebra}
\citation{sagemath}
\@writefile{loa}{\contentsline {algorithm}{\numberline {4}{\ignorespaces Relative Couveignes' method for $e$-th root modulo \(p\)\relax }}{6}{algorithm.4}\protected@file@percent }
\newlabel{alg:crt-couveignes-one-prime}{{4}{6}{Relative Couveignes' method for $e$-th root modulo \(p\)\relax }{algorithm.4}{}}
\newlabel{alg:crt-couveignes-one-prime@cref}{{[algorithm][4][]4}{[1][6][]6}}
\newlabel{lem:gener-couv-one-prime}{{6}{6}{}{lemma.6}{}}
\newlabel{lem:gener-couv-one-prime@cref}{{[lemma][6][]6}{[1][6][]6}}
\@writefile{loa}{\contentsline {algorithm}{\numberline {5}{\ignorespaces Relative Couveignes' method for $e$-th root\relax }}{6}{algorithm.5}\protected@file@percent }
\newlabel{alg:crt-couveignes}{{5}{6}{Relative Couveignes' method for $e$-th root\relax }{algorithm.5}{}}
\newlabel{alg:crt-couveignes@cref}{{[algorithm][5][]5}{[1][6][]6}}
\@writefile{toc}{\contentsline {subsection}{\numberline {5.2}A recursive relative Couveignes' method}{6}{subsection.5.2}\protected@file@percent }
\newlabel{sec:recurs-relat-couv}{{5.2}{6}{A recursive relative Couveignes' method}{subsection.5.2}{}}
\newlabel{sec:recurs-relat-couv@cref}{{[subsection][2][5]5.2}{[1][6][]6}}
\newlabel{prop:recursive-couveignes-applicability}{{4}{6}{}{proposition.4}{}}
\newlabel{prop:recursive-couveignes-applicability@cref}{{[proposition][4][]4}{[1][6][]6}}
\@writefile{toc}{\contentsline {subsection}{\numberline {5.3}Experimental results}{6}{subsection.5.3}\protected@file@percent }
\newlabel{sec:couveignes-simple-exp}{{5.3}{6}{Experimental results}{subsection.5.3}{}}
\newlabel{sec:couveignes-simple-exp@cref}{{[subsection][3][5]5.3}{[1][6][]6}}
\citation{trager}
\citation{fieker-reconstruction,belabas_roots}
\citation{multiquadratics,norm_relations,tw-sti}
\citation{bernard-prepro,tw-sti}
\citation{norm_relations}
\citation{BEFHY22}
\@writefile{lof}{\contentsline {figure}{\numberline {2}{\ignorespaces Timings (s) for \texttt  {nfroots} and\nonbreakingspace Alg.\nonbreakingspace \ref  {alg:crt-couveignes} plotted against \(n\), over fields \(\mathbb  {Q}(\zeta _{pq})\) with constant \([K:L]= p-1\) and \(e=q\).\relax }}{7}{figure.caption.9}\protected@file@percent }
\newlabel{fig:couveignes-naive-reldeg}{{2}{7}{Timings (s) for \texttt {nfroots} and~Alg.~\ref {alg:crt-couveignes} plotted against \(n\), over fields \(\Q (\zeta _{pq})\) with constant \([\nfsup :\nfsub ]= p-1\) and \(e=q\).\relax }{figure.caption.9}{}}
\newlabel{fig:couveignes-naive-reldeg@cref}{{[figure][2][]2}{[1][6][]7}}
\@writefile{lof}{\contentsline {subfigure}{\numberline{(a)}{\ignorespaces {\(p = 5\)}}}{7}{figure.caption.9}\protected@file@percent }
\@writefile{lof}{\contentsline {subfigure}{\numberline{(b)}{\ignorespaces {\(p = 11\)}}}{7}{figure.caption.9}\protected@file@percent }
\@writefile{lof}{\contentsline {subfigure}{\numberline{(c)}{\ignorespaces {\(p = 23\)}}}{7}{figure.caption.9}\protected@file@percent }
\@writefile{lof}{\contentsline {figure}{\numberline {3}{\ignorespaces Timings (s) for \texttt  {nfroots} and\nonbreakingspace Alg.\nonbreakingspace \ref  {alg:crt-couveignes} plotted against \(n\), over fields \(\mathbb  {Q}(\zeta _{pq})\) with constant \(e=p\) and \([K:L]=\varphi (q)\).\relax }}{7}{figure.caption.10}\protected@file@percent }
\newlabel{fig:couveignes-naive-root}{{3}{7}{Timings (s) for \texttt {nfroots} and~Alg.~\ref {alg:crt-couveignes} plotted against \(n\), over fields \(\Q (\zeta _{pq})\) with constant \(e=p\) and \([\nfsup :\nfsub ]=\varphi (q)\).\relax }{figure.caption.10}{}}
\newlabel{fig:couveignes-naive-root@cref}{{[figure][3][]3}{[1][6][]7}}
\@writefile{lof}{\contentsline {subfigure}{\numberline{(a)}{\ignorespaces {\(p = 5\)}}}{7}{figure.caption.10}\protected@file@percent }
\@writefile{lof}{\contentsline {subfigure}{\numberline{(b)}{\ignorespaces {\(p = 11\)}}}{7}{figure.caption.10}\protected@file@percent }
\@writefile{lof}{\contentsline {subfigure}{\numberline{(c)}{\ignorespaces {\(p = 23\)}}}{7}{figure.caption.10}\protected@file@percent }
\@writefile{toc}{\contentsline {section}{\numberline {6}Saturation: a real-life example}{7}{section.6}\protected@file@percent }
\newlabel{sec:saturation}{{6}{7}{Saturation: a real-life example}{section.6}{}}
\newlabel{sec:saturation@cref}{{[section][6][]6}{[1][7][]7}}
\@writefile{toc}{\contentsline {subsection}{\numberline {6.1}Detecting \(e\)-th powers}{7}{subsection.6.1}\protected@file@percent }
\citation{Schiro}
\citation{Schiro}
\citation{number_field_sieve}
\citation{Storjohann}
\citation{tw-sti}
\bibstyle{ACM-Reference-Format}
\bibdata{biblio.bib}
\bibcite{AMM77}{{1}{1977}{{Adleman et~al\mbox  {.}}}{{Adleman, Manders, and Miller}}}
\@writefile{toc}{\contentsline {subsubsection}{\numberline {6.1.1}Conditions on the primes}{8}{subsubsection.6.1.1}\protected@file@percent }
\@writefile{toc}{\contentsline {subsubsection}{\numberline {6.1.2}Definition of the characters}{8}{subsubsection.6.1.2}\protected@file@percent }
\@writefile{toc}{\contentsline {subsubsection}{\numberline {6.1.3}Number of characters}{8}{subsubsection.6.1.3}\protected@file@percent }
\@writefile{toc}{\contentsline {subsection}{\numberline {6.2}Practical considerations}{8}{subsection.6.2}\protected@file@percent }
\@writefile{toc}{\contentsline {subsubsection}{\numberline {6.2.1}Computing suitable primes}{8}{subsubsection.6.2.1}\protected@file@percent }
\@writefile{toc}{\contentsline {subsubsection}{\numberline {6.2.2}Computing non trivial relations}{8}{subsubsection.6.2.2}\protected@file@percent }
\@writefile{toc}{\contentsline {subsubsection}{\numberline {6.2.3}Complexity analysis}{8}{subsubsection.6.2.3}\protected@file@percent }
\@writefile{toc}{\contentsline {subsection}{\numberline {6.3}Experiments}{8}{subsection.6.3}\protected@file@percent }
\bibcite{multiquadratics}{{2}{2017}{{Bauch et~al\mbox  {.}}}{{Bauch, Bernstein, de~Valence, Lange, and van Vredendaal}}}
\bibcite{belabas_roots}{{3}{2004}{{Belabas}}{{Belabas}}}
\bibcite{tw-sti}{{4}{2022}{{Bernard et~al\mbox  {.}}}{{Bernard, Lesavourey, Nguyen, and Roux-Langlois}}}
\bibcite{bernard-prepro}{{5}{2020}{{Bernard and Roux-Langlois}}{{Bernard and Roux-Langlois}}}
\bibcite{BEFHY22}{{6}{2022}{{Biasse et~al\mbox  {.}}}{{Biasse, Erukulangara, Fieker, Hofmann, and Youmans}}}
\bibcite{norm_relations}{{7}{2020}{{Biasse et~al\mbox  {.}}}{{Biasse, Fieker, Hofmann, and Page}}}
\bibcite{Bre75}{{8}{1975}{{Brent}}{{Brent}}}
\bibcite{number_field_sieve}{{9}{1993}{{Buhler et~al\mbox  {.}}}{{Buhler, Lenstra, and Pomerance}}}
\bibcite{cohen-course}{{10}{1993}{{Cohen}}{{Cohen}}}
\bibcite{cohen-advanced}{{11}{2012}{{Cohen}}{{Cohen}}}
\bibcite{couveignes}{{12}{1997}{{Couveignes}}{{Couveignes}}}
\bibcite{fieker-reconstruction}{{13}{2000}{{Fieker and Friedrichs}}{{Fieker and Friedrichs}}}
\bibcite{GN08}{{14}{2008}{{Gama and Nguyen}}{{Gama and Nguyen}}}
\bibcite{GG13}{{15}{2013}{{{Joachim von zur Gathen and J{\"{u}}rgen Gerhard}}}{{{Joachim von zur Gathen and J{\"{u}}rgen Gerhard}}}}
\bibcite{lang2012algebra}{{16}{2012}{{Lang}}{{Lang}}}
\bibcite{Molin}{{17}{2010}{{Molin}}{{Molin}}}
\bibcite{neukirch}{{18}{1999}{{Neukirch}}{{Neukirch}}}
\bibcite{PARI2}{{19}{2022}{{{The PARI~Group}}}{{{The PARI~Group}}}}
\bibcite{sagemath}{{20}{2023}{{{The Sage Developers}}}{{{The Sage Developers}}}}
\bibcite{Schiro}{{21}{1993}{{Schirokauer}}{{Schirokauer}}}
\bibcite{Sho93}{{22}{1993}{{Shoup}}{{Shoup}}}
\bibcite{Storjohann}{{23}{2013}{{Storjohann}}{{Storjohann}}}
\bibcite{thome}{{24}{2012}{{Thom{\'e}}}{{Thom{\'e}}}}
\bibcite{trager}{{25}{1976}{{Trager}}{{Trager}}}
\newlabel{tocindent-1}{0pt}
\newlabel{tocindent0}{0pt}
\newlabel{tocindent1}{4.185pt}
\newlabel{tocindent2}{10.34999pt}
\newlabel{tocindent3}{18.198pt}
\@writefile{lof}{\contentsline {figure}{\numberline {4}{\ignorespaces Timings (s) for \texttt  {nfroots} and Alg. 4 plotted against the dimension for saturation process.\relax }}{9}{figure.caption.11}\protected@file@percent }
\newlabel{fig:couveignes-saturation}{{4}{9}{Timings (s) for \texttt {nfroots} and Alg. 4 plotted against the dimension for saturation process.\relax }{figure.caption.11}{}}
\newlabel{fig:couveignes-saturation@cref}{{[figure][4][]4}{[1][8][]9}}
\@writefile{lof}{\contentsline {subfigure}{\numberline{(a)}{\ignorespaces {\(e \not \mid m\)}}}{9}{figure.caption.11}\protected@file@percent }
\@writefile{lof}{\contentsline {subfigure}{\numberline{(b)}{\ignorespaces {\(e \mid m\)}}}{9}{figure.caption.11}\protected@file@percent }
\@writefile{toc}{\contentsline {section}{References}{9}{section*.13}\protected@file@percent }
\newlabel{TotPages}{{9}{9}{}{page.9}{}}
\gdef \@PBprevLastPage{9}
\csdimgdef {@PBunbalpgIXHeight}{481.20726pt}
\csdimgdef {@PBunbalpgIXLeftHeight}{481.20726pt}
\csdimgdef {@PBunbalpgIXLeftFloatsHeight}{0.0pt}
\csdimgdef {@PBunbalpgIXRightHeight}{481.20726pt}
\csdimgdef {@PBunbalpgIXRightFloatsHeight}{0.0pt}
\csdimgdef {@PBunbalpgIXUsedLeft}{476.75pt}
\csdimgdef {@PBunbalpgIXUsedRight}{34.0pt}
\providetoggle {@PBunbalpgIXHasFloats}
\global \toggletrue {@PBunbalpgIXHasFloats}
\providetoggle {@PBunbalpgIXHasFloatcol}
\global \togglefalse {@PBunbalpgIXHasFloatcol}
\providetoggle {@PBunbalpgIXHasFootnotes}
\global \togglefalse {@PBunbalpgIXHasFootnotes}
\providetoggle {@PBunbalpgIXHasMarginpars}
\global \togglefalse {@PBunbalpgIXHasMarginpars}
\ExplSyntaxOn
\seq_gclear_new:c{@PBunbalpgIXLeftFloatHeights}
\ExplSyntaxOff
\ExplSyntaxOn
\seq_gclear_new:c{@PBunbalpgIXLeftFloatSpacesBelow}
\ExplSyntaxOff
\providetoggle {@PBpgIXalreadyBalanced}
\global \toggletrue {@PBpgIXalreadyBalanced}
\csdimgdef {@PBbalpgIXHeight}{481.20726pt}
\csdimgdef {@PBbalpgIXLeftHeight}{481.20726pt}
\csdimgdef {@PBbalpgIXLeftFloatsHeight}{0.0pt}
\csdimgdef {@PBbalpgIXRightHeight}{481.20726pt}
\csdimgdef {@PBbalpgIXRightFloatsHeight}{0.0pt}
\csdimgdef {@PBbalpgIXUsedLeft}{267.375pt}
\csdimgdef {@PBbalpgIXUsedRight}{258.0pt}
\providetoggle {@PBbalpgIXHasFloats}
\global \toggletrue {@PBbalpgIXHasFloats}
\providetoggle {@PBbalpgIXHasFloatcol}
\global \togglefalse {@PBbalpgIXHasFloatcol}
\providetoggle {@PBbalpgIXHasFootnotes}
\global \togglefalse {@PBbalpgIXHasFootnotes}
\providetoggle {@PBbalpgIXHasMarginpars}
\global \togglefalse {@PBbalpgIXHasMarginpars}
\ExplSyntaxOn
\seq_gclear_new:c{@PBbalpgIXLeftFloatHeights}
\ExplSyntaxOff
\ExplSyntaxOn
\seq_gclear_new:c{@PBbalpgIXLeftFloatSpacesBelow}
\ExplSyntaxOff
\providetoggle {@PBstabilized}
\global \toggletrue {@PBstabilized}
\providetoggle {@PBimpossible}
\global \togglefalse {@PBimpossible}
\gdef \@PBprevLastPhysPage{9}
\gdef \@abspage@last{9}
\makeatother}{}%
}
\begin{document}

\begin{abstract}
  We describe several algorithms for computing $e$-th roots of elements in a number
  field $K$, where $e$ is an odd prime integer. In particular we generalize
  Couveignes' and Thom\'e's algorithms originally designed to compute square-roots
  in the Number Field Sieve algorithm for integer factorization. Our algorithms cover 
  most cases of $e$ and $K$ and allow to obtain reasonable timings even for large 
  degree number fields and large exponents $e$. The complexity of our algorithms is 
  better than general root finding algorithms and our implementation compared well 
  in performance to these algorithms implemented in well-known computer algebra softwares. One 
  important application of our algorithms is to compute the saturation phase in the 
  Twisted-PHS algorithm for computing the Ideal-SVP problem over cyclotomic fields 
  in post-quantum cryptography. 
\end{abstract}

\maketitle
\addtolength{\leftmargini}{-1ex}

\section{Introduction}
\label{sec:introduction}
Computing roots of elements is an important step when solving various tasks in
computational number theory. It arises for example during the final step of the
General Number Field Sieve~\cite{number_field_sieve,couveignes} (NFS). This
problem also intervenes during saturation processes while computing the class
group or $S$-units of a number field~\cite{norm_relations}. Recently,
such computations were found to be important to study the Ideal-SVP problem used
in Lattice-based Cryptography~\cite{tw-sti}.

Generally speaking, elements are given in a ``factored'' form, meaning that an
element $y$ for which a root needs to be computed is given as a product of
relatively small elements $y = \prod_{i=1}^r u_i^{e_i}$. Note that the two
contexts mentioned previously are somewhat orthogonal ones to each other. Indeed,
for the NFS, the length $r$ of the product is very large while the degree of the
number field is typically small, about 10, and one needs to compute square roots
($e = 2$). In saturation processes, such as the ones we are interested in
(see~\cite{tw-sti} for practical cases), $r$ is typically a few times the degree
of the field, which is potentially large, say between 100 and 200. Moreover the
exponent \(e\) may be very large as well, $90$ bits. Thus, most strategies to
compute an \(e\)-th root developed in the NFS context 
become intractable in this setting if not carefully adapted.

\subsection*{Our contributions}

\def\nfsup{K}
\def\nfsub{L}
\def\conddblcrt{($*$)}

In this article, we explain how to efficiently compute an $e$-th root of
an element $y \in K$, where $K$ is a number field and \(e\) is an odd prime
power.
We aim at designing a workable method for \emph{large} exponents~$e$ and
dimension $[K:\Q]$ for all cases.
\begin{itemize}
\item When $K$ and $e$ are such that there exist infinitely many prime integers $p$ such that
  $\forall \mf p \mid p$, $p^{f_{\mf p}} \not \equiv 1 \bmod e$ (condition \conddblcrt), 
  we reconstruct \(x\) from \( \bigl(x \bmod {p_1}, \dotsc, x \bmod {p_k}\bigr) \) using the
  Chinese Remainder Theorem (CRT), where each \((x \bmod{p_i})\) is itself
  computed through a CRT procedure on prime ideals of $K$. We call
  this generalisation of Thom\'e's square-roots algorithm~\cite{thome}
  the Double-\CRT{}~(\cref{alg:crt-primes,alg:crt-root}).
\item Generically, we can use $p$-adic lifting when there is an inert prime, or
  the $\mf{p}$-adic reconstruction of
  Belabas~\cite{belabas_roots}.
  Both use a
  Hensel's lifting, that we adapt to compute $e$-th roots while avoiding inverse
  computations in \S\ref{sec:generic} (see \cref{alg:padic}). While both methods work for any $e$,
  inert primes do not~always exist and $\mf p$-adic reconstruction scales poorly
  with the dimension of $K$.
  However, sometimes they can be very useful and we indeed use them in our
  last recursive algorithm.
\item When good conditions on \(K\) and \(e\) are not satisfied, we show how one
  can adapt Couveignes' approach for square roots~\cite{couveignes} to relative
  extensions of number fields \(\nfsup/\nfsub\), provided \([\nfsup:\nfsub]\) is
  coprime to \(e\) and sufficiently many prime integers \(p\) verify that each 
  prime ideal \(\mf p\) of \(\O \nfsub\) above \(p\) is inert in \(\nfsup\). We then
  turn this strategy into a recursive procedure, calling the previous algorithms
  until the smallest possible subfield is reached~(\cref{alg:crt-couveignes-one-prime,alg:crt-couveignes}).
\end{itemize}
This leads to the following global strategy to compute an \(e\)-th root of an
element \(y \in \nfsup\), which we use in our implementation:
\begin{enumerate}
\item if condition \conddblcrt{} is satisfied,
  use the Double-CRT (\cref{alg:crt-root});
\item else, if inert primes exist in $K$, use a \(p\)-adic lifting
  (see \S\ref{padiclifting});
\item otherwise, try to use the relative Couveignes'
  strategy, going back to step\,(2) for computing an~\(e\)-th root of
  \(\norm{\nfsup/\nfsub}(y)\) in \(\nfsub\);
\item finally, resort to the \(\mf p\)-adic reconstruction given in \S\ref{padicreconstruction}.
\end{enumerate}

\subsection*{Experimental results}
We ran experiments to evaluate the performances of our algorithms when compared
to standard methods and implementations, especially \textsc{Pari/Gp}
\textsf{nfroots}. We focused on cyclotomic fields, as they are the main fields
used in our application domain, \ie lattice-based cryptography, but our
algorithms extend to other number fields in most cases. All of our
implementations are done using \textsc{SageMath}~\cite{sagemath} with few
optimisations. Meanwhile, we compare with \textsc{Pari/Gp}, and still achieve
several orders of performance (between 10 to 10000+ when $n$ and $e$
increase). Thus, our timings could be further improved with an optimised C
implementations.
All our codes are freely available at \url{\giturl}.

Over ``good'' cases, our~\CRT; generalisation algorithm is clearly more
efficient than \textsc{Pari/Gp} \textsf{nfroots},
see e.g.~\Cref{fig:thome-naive-root}, and the gap explodes when the exponent
\(e\) increases. 
Our algorithm scale well
and we can use it without any problems for $e$ of $94$-bits prime (see~\cref{tab:timings-s-units}),
which would be completely irrealistic for \textsc{Pari/Gp}.

Over ``bad'' cases, experiments show that our generalization of Couveignes'
algorithm is also more efficient than \textsc{Pari/Gp} \textsf{nfroots}~\cite{PARI2},
see \Cref{fig:couveignes-naive-reldeg,fig:couveignes-naive-root}. Our algorithm
is always faster, and the gap with \textsf{nfroots} becomes larger when $e$ and
$n$ increase.

Finally, we tested our algorithms in a real-life situation, namely saturating
full-rank multiplicative sets of $S$-units arising from Stickelberger’s theorem
using the code of~\cite{tw-sti}, see \S\ref{sec:saturation},
\cref{tab:timings-s-units,fig:couveignes-saturation}. Again, our implementation
of our algorithms
is more efficient than \textsc{Pari/Gp} \textsf{nfroots} for all ranges of
exponents or dimensions.


\section{Preliminaries}
\label{sec:preliminaries}


Let $K = \Q(\alpha)$ be a number field defined by a monic irreducible polynomial
$f\in\Z[t]$ of degree $n$ s.t.~$f(\alpha)=0$.  In this paper, we shall suppose
that we know a reasonable (e.g.,~LLL-reduced) basis $(\omega_i)$ of some order
$\Z[\alpha]\subseteq\O{}\subseteq\O{K}$ and $f_0\geq 1$ such that
$f_0\O{K}\subset\O{}$.

\paragraph{Complexities.} We use the standard $O(\cdot)$ notation. Let us denote by
$\Mint{s}$ the complexity of multiplying two $s$-word integers,
and by $\Mpol{d}{s}$ the complexity of multiplying two polynomials of degree $d$ whose
coefficients are taken modulo an $s$-word integer, \ie using fast arithmetic (see \eg \cite[Cor.\,11.8~and\,11.10]{GG13}),
\begin{equation*}
  \Mpol{d}{s} \leq O(d\log d \llog d \cdot s \log^2 s \llog s).
\end{equation*}

\subsection{Bounds on root coefficients}
\label{sec:bound}



In all the methods of this paper, we need a bound on the coefficients of the
seeked $e$-th root on the basis $(\omega_i)$.
Such bounds are usually obtained by computing a lazy estimation of the inverse
of a Vandermonde-like matrix linking complex embeddings of elements of $K$ to
their coefficients corresponding to $(\omega_i)$.

More precisally, let $\Omega = (b_{ij}) \in \mathcal{M}(\Q)$ be
s.t.~$\omega_i = \sum_j b_{ij}\alpha^{j-1}$, and
$V_{\alpha} = \bigl( \sigma_j(\alpha)^{i-1}\bigr)$ the Vandermonde matrix
corresponding to $f$. For $x = \sum_i c_i \omega_i$, let $C(x) = (c_i)$ be the
coefficient embedding of $x$ and
$\Sigma(x)= C(x) \cdot \bigl(\Omega V_{\alpha}\bigr)$ be its canonical embedding
to $\CC$.
\begin{lemma}\label{lem:boundcoef}
  Define by $\Ninf{A} = \max_{j} \sum_i |a_{ij}|$ the infinity norm of a matrix $A$, and let $\Cinf = \Ninf{V_{\alpha}^{-1} \Omega^{-1}}$. Then we have
\begin{equation*}
  \Ninf{C(x)} \leq \Ninf{\Sigma(x)} \cdot \Cinf.
\end{equation*}
\end{lemma}

\begin{proof}
  This is a direct adaptation of \cite[Lem.\,3.3]{belabas_roots}
  (see also \cite[Lem.\,6]{fieker-reconstruction}), noting that
  the given definition of the infinity norm for matrices is equivalent to
  $\Ninf{A} = \sup_{y\neq 0} \bigl( \Ninf{yA} / \Ninf{y} \bigr)$.
\end{proof}


Hence, for any $x\in K$, the coefficient norm $\Norm{C(x)}$ is only a constant
factor away from the usual norm $\Norm{\Sigma(x)}$, and in practice only a loose
estimation of this factor is necessary.
In particular, for~$y\in (K^*)^e$, it is easy to evaluate the size of the
(canonical) embedding norm of its $e$-th root $x$ as
$\ln \Ninf{\Sigma(x)} = \tfrac 1e \ln \Ninf{\Sigma(y)}$, and
\cref{lem:boundcoef} gives the intuition that generically, the coefficients
of~$x$ are roughly~$e$ times smaller than those of $y$, which holds well in
practice.

Moreover, in our particular case of interest, $y$ is given in \emph{factored}
form, \ie there exist $(u_i)_{1\leq i\leq r}\in K^*$ and
$(a_i)_{1\leq i\leq r} \in \segment{0}{e-1} ^ r$ s.t.~
 $ y = \textstyle\prod_{1\leq i \leq r} u_i^{a_i}$.
 In this situation, the following lemma shows that the size of the coefficients
 of $x = y^{1/e}$ does \emph{not} depend on $e$, but only on the total size of
 the $u_i$'s. 
\begin{lemma}\label{lem:boundethroot}
  Let $y = \prod_i u_i^{a_i}$ be an element of $(K^*)^e$ in factored form as
  above, and let $x$ be such that $y=x^e$. Then
  \begin{equation*}
    \ln \Ninf{\Sigma(x)} < \textstyle\sum_{1\leq i\leq r} \ln \Ninf{\Sigma(u_i)}. 
  \end{equation*}
\end{lemma}
\begin{proof}
  For any embedding $\sigma$, we have
  \begin{equation*}
    \ln \abs{\sigma(x)} = \tfrac{1}{e} \cdot \ln \abs{\sigma(y)} = \tfrac{1}{e} \cdot \textstyle\sum_{1\leq i \leq r} a_i \ln \abs{\sigma(u_i)}.
  \end{equation*}
  As $a_i < e$ for all $1\leq i \leq r$, the lemma follows.
\end{proof}

This motivates the design of algorithms that \emph{never} compute $y$ globally,
which is absolutely crucial when $e$ is very large.

\subsection{Computing \texorpdfstring{$e$}{e}-th roots in finite fields}
\label{ethfp}



Let $\F_q$ be a finite field of characteristic $p>2$ with $q=p^d$. For the
local-global methods of this paper, an important tool consists in computing an $e$-th root in
$\F_q$, where $e$ is given modulo $q-1$.

Let $y = x^e$ in $\F_q^*$. 
If $q \not\equiv 1 \bmod{e}$, then every element $y$ is an $e$-th power and we
can simply compute $x$ as $\smash[t]{y^{e^{-1}\bmod (q-1)}}$.
If $q \equiv 1 \bmod{e}$, then we need to work in the subgroup of order $(q-1)/e$ of $\smash{\F_q^*}$.
If~$q\not\equiv 1\mod{e^2}$, we can again compute $x$ as $y^{e^{-1}\bmod
  (q-1)/e}$. Both cases can be treated efficiently in
$O\bigl(\log q \cdot \Mpol{d}{\log p}\bigr)$ operations.

In very rare cases, we will have no choice but to consider
$q \equiv 1 \mod e^2$, thus we resort to the Adleman-Manders-Miller algorithm
\cite[Th.\,IV]{AMM77}, which is an adaptation of Tonnelli-Shanks algorithm to
the case $e>2$, and has complexity
$O\bigl( (e \log q + \log^2 q) \cdot \Mpol{d}{\log p} \bigr)$ in the worst case.
%
Alternatively, we can use generic factorisation methods such as the
Cantor-Zassenhaus algorithm \cite[Alg.\,3.4.6]{cohen-course}, whose complexity
in $O\bigl(e^2 \log^{1+\epsilon} e \log (qe) \cdot \Mpol{d}{\log p}\bigr)$
\cite{Sho93} might be competitive for small $e < \log q$.

For simplicity's sake, we will denote by $\FFroot{e}{q}$ the complexity of computing an $e$-th root in $\F_{q}$ in all cases.




\section{Generic local-global methods}
\label{sec:generic}
In this section, we present two local-global methods that apply unconditionnally \wrt $e$, namely $p$-adic lifting (\eg \cite[\S1.1]{thome}), and $\mf{p}$-adic reconstruction \cite{belabas_roots}.
Both methods rely on a $p$-adic (\resp $\mf{p}$-adic) Newton iteration, or Hensel's lift, that we tweak in the particular case of $e$-th roots to avoid $p$-adic (\resp $\mf{p}$-adic) inversions.

Despite their genericity, both methods have severe drawbacks.
The $p$-adic lifting relies on the existence of inert primes, which rules out many families of number fields, \eg cyclotomic fields of composite conductors.
The $\mf{p}$-adic reconstruction does not scale well as the degree of the number field grows, since it requires LLL-reducing a possibly badly-skewed ideal lattice.

Nevertheless, the techniques developed here will be used as a base case for our (recursive) relative Couveignes' method \S\ref{sec:recurs-relat-couv}, which aims at reducing the dimension of the $e$-th root computation.




\subsection{\texorpdfstring{$p$}{p}-adic lifting}
\label{padiclifting}




The classical $p$-adic lifting approach applies whenever there exist inert prime ideals $p$ in $K$. 
%
It relies on the fact that $K_p$, the $p$-adic completion~of~$K$ at $p$, is then
a degree~$n=[K:\Q]$ unramified extension of $\Q_p$, and the morphism sending a
root in $K$ of the defining polynomial~$f$ to a root of~$f$ in~$K_p$ is
injective. Thus, a sufficiently good approximation in~$K_p$ of the image of
$x\in K$ allows us to retrieve directly the coefficients of~$x$. This is done by
means of Newton iterations in $K_p$.

Let $h(z)= z^e - y$, with $y = x^e$ for some $x\in K^*$, and let $p$ be an inert prime of $K$ with $\gcd(p,e)=1$.
The Newton iteration for $h$ writes as follows.
%
At each step $i\geq 0$, suppose that we know a $p$-adic approximation $x_i$ of
$x$ at precision $k=2^i$, \ie\mbox{$x_i = x + O(p^k)$}. For~$i=0$, this
approximation is found by usual methods in $\smash[b]{\F_{p^n}}$
(\S\ref{ethfp}). Then, $h(x_i) = O(p^k)$ and an approximation $x_{i+1}$ at
precision $2k$ is obtained by computing the following iteration modulo $p^{2k}$:
\begin{equation*}
  x_{i+1} = x_i - \tfrac1e\cdot \Bigl(x_i - \tfrac y{x_i^{e-1}} \Bigr).
\end{equation*}
Let $B>0$ be an upper bound on the coefficients of $x$, then this iteration is
performed $\kappa = \bigl\lceil \log_2 \max \{1, \log_p 2B\}\bigr\rceil$
times. 
In particular, only the knowledge of $y$ at precision $p^{\smash[t]{2^\kappa}}$
is needed. This is especially useful when $e$ is large, since in that case the
size of the coefficients of $x$ is roughly $e$ times smaller than for $y$ (see
\S\ref{sec:bound}).

In practice, each iteration above requires an inverse computation. For $e=2$, a
known trick to avoid this consists in first computing the \emph{inverse} square
root, then using $y^{1/2} = y\cdot y^{-1/2}$ \cite[\S2]{Bre75}. We adapt this
trick in our case by seeking first a root $x$ of the polynomial
\begin{equation*}
  g(z) = y^{e-1} z^e - 1.
\end{equation*}
Then $(y\cdot x)$ verifies $(yx)^e = y\cdot (y^{e-1}x^e) = y$ as expected. The Newton iteration for $g$ now writes without inverses as
\begin{equation}\label{newtiteration}
  x_{i+1} = x_{i} - \tfrac 1e x_i \bigl( y^{e-1} x_i^e - 1 \bigr).
\end{equation}
The complete $p$-adic lifting using this new iteration is summarized in
\cref{alg:padic}.

\begin{algorithm}[H]
  \caption{$p$-adic lifting for $e$-th root}
  \label{alg:padic}
  \begin{algorithmic}[1]
    \Require $y\in (K^*)^e$ in factored form $y = \prod_i u_i^{a_i}$, $p$ an inert prime in $K$ with $\gcd(p,e)=1$.
    \Ensure  \(x\in K^*\) such that $y = x^e$.
    \State Compute $B$ s.t.~$\Ninf{C(x)}\leq B$ \Comment{Using \cref{lem:boundcoef,lem:boundethroot}}
    \State $\kappa \gets \lceil \log_2 \max\{ 1, \log_p 2B\}\rceil$
    \State $a \gets \bigl(\prod_i u_i^{a_i}\bigr)^{e-1} \mod p^{2^{\kappa}}$ \Comment{Reduce $u_i$'s first}
    \State $x_0 \gets (a \mod p)^{1/e}$ in $\F_{p^n} \simeq \faktor{\O{}}{p\O{}}$ \Comment{Using \S\ref{ethfp}}
    \For{$0\leq i < \kappa$}
    \State $x_{i+1} \gets x_{i} - \tfrac 1e x_{i} \bigl(a x_{i}^e - 1\bigr)$ \Comment{Work in $\O{}$ modulo $p^{2^{i+1}}$}
    \EndFor
    \State $x \gets a\cdot x_{\kappa} \mod p^{2^\kappa}$
    \State \Return $x$ with coefficients mapped in $[-B, B]$
  \end{algorithmic}
\end{algorithm}




\begin{proposition}
  \Cref{alg:padic} is correct, and runs in time at most
  $O\bigl(\log e \cdot (r + \log s) \cdot \Mpol{n}{s} + \FFroot{e}{p^n} \bigr)$, where
  $s$ is the total input size, \ie
  $s = O\bigl(\sum_i \log \Ninf{\Sigma(u_i)}\bigr)$.
\end{proposition}

\begin{proof}
  Let $x$ be a root of $g(z) = y^{e-1} z^e - 1$. At any stage $i\geq 0$, let $x_i$ be a $p$-adic approximation of $x$ at precision $k=2^i$, \ie $x = x_i + O(p^k)$, and let $\varepsilon_i=g(x_i) = O(p^k)$. We shall show that $\varepsilon_{i+1} = g(x_{i+1}) = O(p^{2k})$, which implies correctness. Working modulo $p^{2k}$, and plugging the Newton iteration formula for $x_{i+1}$ into $g$, we get
  \begin{equation*}
    \begin{split}
      \varepsilon_{i+1} & = y^{e-1} \bigl( x_{i} - \tfrac 1e x_i g(x_i)\bigr)^e - 1 
                          = y^{e-1} \bigl( x_{i} - \tfrac 1e x_i \varepsilon_i\bigr)^e -1 \\
                        & = - 1 + y^{e-1} x_i^e \bigl( 1 - e \tfrac 1e \varepsilon_i + O(\varepsilon_{i}^2)\bigr)\\
                        & = - 1 + (\varepsilon_i+1) - \varepsilon_i(\varepsilon_i+1) + O\bigl(\varepsilon_{i}^2\bigr) 
                          = O\bigl(\varepsilon_i^2\bigr) = O\bigl(p^{2k}\bigr).
    \end{split}    
  \end{equation*}
  As for the complexity, note that, by \cref{lem:boundethroot}, $\log B = O(s)$, thus \mbox{$\kappa = O(\log s)$}.
  Computing the degree $n$ polynomial $(a \bmod p^{2^{\kappa}})$  
  costs at most $O\bigl( r \log e \cdot \Mpol{n}{s}\bigr)$ using that $a_i < e$, and all polynomials $(a \bmod p^{2^i})$ for the loop can be iteratively deduced in negligible $O(n \Mint{s}\log s)$ time. Likewise, computing all $(1/e \bmod p^{2^i})$ costs $O(s \log e \log s)$.
  Computing $(a\bmod p)^{1/e}$ in $\F_{p^n}$ costs $\FFroot{e}{p^n}$ by \S\ref{ethfp}.
  Finally, using these values, each of the $O(\log s)$ Newton iterations at precision $k=2^i$ costs at most $O\bigl(\log e \cdot\Mpol{n}{k\log p}\bigr)$ for a total of $O\bigl(\log e \log s \cdot\Mpol{n}{s}\bigr)$.
\end{proof}







\subsection{\texorpdfstring{$\mf{p}$}{P}-adic reconstruction}
\label{padicreconstruction}
When the field $K$ contains no inert primes, the above method can still be used,
for any unramified prime ideal $\mf{p}$ of inertia degree $\rd{\mf{p}}{p} < n$,
to obtain a $\mf{p}$-adic approximation in $K_{\mf{p}}$, the completion of $K$
at $\mf{p}$, which is an unramified extension of $\Q_p$ of degree
$\rd{\mf{p}}{p}$.

Starting from a low-precision $e$-th root in
\(\smash{\faktor{\O{K}}{\mf{p}} \simeq
  \F_{p^{\rd{\mf{p}}{p}}}}\),~\cref{newtiteration} allows for its lifting modulo
$\mf{p}^{a}$ for any $a$.  If $a$ is large enough, it is possible to reconstruct
the root in $K$ from this $\mf{p}$-adic embedding approximation, as is done in
\cite[\S3]{belabas_roots}, by solving a Bounded Distance Decoding problem in the
LLL-reduced lattice corresponding to $\mf{p}^a$. The main drawback of this
method is that the ideal $\mf{p}^a$ is all the more badly-skewed that the
inertia degree of $\mf{p}$ is small and $n$ is large, so that the LLL-reduction
quickly dominates \cite[\S3.7]{belabas_roots}.

In practice, we estimate $a$ using the analysis of \cite{belabas_roots}. Suppose
that we have a bound $B'$ on the coefficient norm of~$x$. By
\cite[Lem.\,3.7]{belabas_roots}, the reconstruction succeeds when
$r_{\max} > B'$, where $r_{\max}$ is explicitly given by
\cite[Lem.\,3.8]{belabas_roots}, which provides a lower bound on $a$ as
in \cite[Lem.\,3.12]{belabas_roots}. Hence, using $\gamma \approx 1.022$ as the
root-Hermite factor achieved by LLL \cite{GN08}, we start from the smallest
$a = 2^\kappa$ such that
\begin{equation*}
  a > \frac{n}{\rd{\mf{p}}{p}\cdot\ln p} \Bigl( \ln 2B' + 
  \bigl(\tfrac{n(n-1)}{4}-1\bigr) \ln \gamma \Bigr).
\end{equation*}
This value of $a$ is controlled \textit{a posteriori} by computing the
corresponding $r_{\max}$, checking whether $r_{\max} > B'$ and doubling $a$
while necessary. In practice, the above estimation is rarely invalidated.


\section{Using Chinese Remainder Theorem: the easy cases}
\label{sec:roots-easy}

In this section, we describe how one can compute $e$-th roots using the Chinese
Remainder Theorem in number fields when $e$
verifies
\begin{equation}
  \label{eq:goodcase}
  \exists_{\infty} \text{ primes $q$ s.t.~}\forall \mf{q}\mid q, q^{f_{\mf{q}}}
  \not\equiv 1 \mod e.
\end{equation}
This condition ensures that none of the residue fields contain a primitive
$e$-th root of unity.  In cyclotomic fields of conductor $m$, this condition is
equivalent to assuming $\gcd(e, m) = 1$, since in that case all primes verify
the assumption and conversely.

In the context of the Number Field Sieve, Thomé described a CRT-based method to
compute square-roots~\cite[\S4]{thome}. A major problem for $e=2$ is to guess
the correct signs modulo each prime ideal, which is handled by solving a
knapsack problem. However, as $n$ grows as well as $e$, this approach quickly
becomes intractable.

\subsection{A CRT-based method for \texorpdfstring{$e$}{e}-th roots}


First, we show in \cref{alg:crt-primes} how to retrieve an $e$-th root modulo $q$, where \(q\) is a prime
integer verifying \eqref{eq:goodcase}, using the Chinese Remainder Theorem in number fields.

\begin{algorithm}[ht]
  \caption{Number field CRT for $e$-th root mod $q$}
  \label{alg:crt-primes}
  \begin{algorithmic}[1]
    \Require An unramified prime $q$ in $K$ verifying \eqref{eq:goodcase}, and
    $y\in (K^*)^e$ in factored form $y = \prod_{i \leqslant r} u_i^{a_i}$, where the $u_i$'s are given mod $q$.
    \Ensure \(x \equiv y^{1/e} \bmod (q)\).
    \State \(S \gets \lbrace \mf{q}_1, \dots, \mf q_g \rbrace := \lbrace \mf{q};\ \mf{q}
           \mid q \O{} \rbrace \) \Comment{Using Cantor-Zassenhaus}
    \State Compute all $(u_i \mod \mf q_j)_{i,j}$ \Comment{Use a product tree}
    \For{\(\mf{q} \in S\)} 
        \State $y_{\mf q} \gets \prod_i (u_i \bmod \mf{q})^{a_i}$ 
        \State $x_{\mf q} \gets {y_{\mf q}^{1/e}} \bmod \mf{q}$ \label{nfcrt:stepmodid} \Comment{Using \S\ref{ethfp}}
    \EndFor
    \State \Return
    \( \text{CRT}_K\bigl(\{x_{\mf{q}} \}_{\mf q \in S}, \{\mf q\}_{\mf q \in S}\bigr)\) \label{nfcrt:idcrt}
  \end{algorithmic}
\end{algorithm}
\begin{algorithm}[htb]
  \caption{Double CRT for $e$-th root}
  \label{alg:crt-root}
  \begin{algorithmic}[1]
    \Require $y\in (K^*)^e$ in factored form $y = \prod_{i \leqslant r} u_i^{a_i}$.
    \Ensure  \(x \in K^*\) such that $y = x^e$.
    \State Compute $B$ s.t.~$\Ninf{C(x)}\leq B$ \Comment{Using \cref{lem:boundcoef,lem:boundethroot}}
    \State Choose primes $q_1, \dotsc, q_k$ verifying \eqref{eq:goodcase} s.t.~$\prod_j q_j \geq 2B$.
    \State Compute all $(u_i \mod q_j)_{i,j}$ \Comment{Use a product tree}
    \For{$1\leq j \leq k$}
        \State $x_j \gets (y \mod q_j)^{1/e}$ 
        \label{dblcrt:stepmodq}
        \Comment{Using \cref{alg:crt-primes}}
    \EndFor
    \State $x \gets \text{CRT}_{\Z}\bigl(\{x_j\}_{j}, \{q_j\}_{j}\bigr)$ %
    \label{dblcrt:lastcrt}%
    \Comment{Coefficient by coefficient}
    \State \Return $x$ with coefficients mapped in $[-B, B]$
  \end{algorithmic}
\end{algorithm}

\begin{proposition}\label{pr:cplxcrtq}
  \Cref{alg:crt-primes} is correct, and runs in time at most
  $O\bigl( r n\log q \cdot (\Mpol{\max \{ f_{\mf q} \}}{\log q} + \log e)\bigr)$.
\end{proposition}

\begin{proof}
  Since~\(q\) is unramified in \(K\), \ie 
  \((q) = \prod_{\mf q \in S} \mf q\), one has 
  \[
    \faktor{\O K}{(q)} \cong \prod_{{\mf q \mid (q)}} \faktor{\O K}{\mf q} \cong
    \prod_{{\mf q \mid (q)}} \F_{q^{\rd{\mf q}{q}}}.
  \]
  For all \(\mf q \in S\), the condition
  \(\smash[t]{q ^ {f_{\mf q}}} \not \equiv 1 \bmod e \) implies
  \(\F_{q^{\smash{f_{\mf q}}}}^* = (\F_{q^{\smash{f_{\mf q}}}}^*)^e\), so that any
  element of \hspace{0.2em}\(\smash{\faktor{\O{}}{\mf q}}\) has a \emph{unique}
  \(e\)-th root.  Consequently, step\,\ref{nfcrt:stepmodid} of \cref{alg:crt-primes} is properly
  defined and \(x \equiv x_{\mf q} \bmod \mf q \). Thus the output of step\,\ref{nfcrt:idcrt} is
  indeed congruent to \(x\) modulo \(q\).

  As for the complexity, the first step, factorising a degree $n$ polynomial
  over \(\F_q\), can be done in
  \(O\bigl(n^2 \log^{1+\epsilon} n \log qn \cdot \Mint{\log q} \bigr)\) using
  Cantor-Zassenhaus.
  Computing the product tree for the $\mf q_j$'s and reducing $\{u_i\}_{1\leq i\leq r}$ modulo all $\mf q_j$'s cost
  $(r+1)\cdot O\bigl(\log g \cdot \Mpol{n}{\log q}\bigr)$ \cite[Lem.\,10.4~and~Th.\,10.15]{GG13}.
  For a given $\mf q \mid q$, computing $y_{\mf q}$ costs at
  most $O\bigl(r f_{\mf q}\log q \cdot \Mpol{f_{\mf q}}{\log q}\bigr)$ (mapping
  each $a_i$ modulo $q^{f_{\mf q}}-1$), and the $e$-th root in
  \( \F_{\smash{q^{f_{\mf q}}}}\) costs
  \( O\bigl(f_{\mf q}\log q\cdot \Mpol{f_{\mf q}}{\log q}\bigr) \) since by hypothesis $q^{f_{\mf q}}\not\equiv 1\bmod e$.
  Since all $f_{\mf q}$ sum to $n$, the whole loop costs at most
  $O\bigl(r n\log q\cdot \Mpol{\max f_{\mf q}}{\log q}\bigr)$.
  Note that if $e$ is large, we can reduce it
  as well as the $a_i$'s once modulo $q^n-1$, and use these smaller versions
  thereafter, yielding an extra $O(r \log e \cdot n \log q)$ term in the worst scenario.
  The last CRT step can be done in $O\bigl(\log n \cdot \Mpol{n}{\log q}\bigr)$
  \cite[Cor.\,10.23]{GG13}, which is negligible.
\end{proof}

\begin{remark}\label{rem:crt-primes}
  Note that the complexity of \cref{alg:crt-primes} does not depend heavily on
  $e$.  Furthermore, the proof shows that it is best to consider primes with
  small maximum inertia degrees. If there are sufficiently many totally split primes,
  \eg in cyclotomic fields, the complexity even drops to 
  $O\bigl(r n \log^{2+\epsilon} q + r \log q \log e\bigr)$.
\end{remark}


We now explain how to compute \(x\) \textit{via} computations modulo primes \(q_1, \dots, q_k\)
that all verify \eqref{eq:goodcase}. This is done in \cref{alg:crt-root}.

\begin{proposition}
  \Cref{alg:crt-root} is correct, and runs in time at most
  $O\bigl( r n  s \cdot ( \Mpol{n}{\kappa} + \log^{2+\epsilon} s + \log e)\bigr)$,
  where 
  $s = O\bigl(\sum_i \log \Ninf{\Sigma(u_i)}\bigr)$
  is the total input size.
\end{proposition}

\begin{proof}
  Step\,\ref{dblcrt:stepmodq} computes the solution modulo the ideal generated by
  \(q_j\), for each \(j \in \segment{1}{k}\). Thus the CRT in step\,\ref{dblcrt:lastcrt} computes
  \(x \mod (Q)\), where $Q = \prod_j q_j$, and \(Q > 2B \) ensures that there is a unique
  element \(z \equiv x \bmod (Q)\) with coefficients in $[-B, B]$.

  Let $\kappa = \max_j \log q_j$.
  To fix ideas, $\kappa$ is taken as to fit a machine word
  and all $q_j$'s are chosen evenly,
  so %
  $k \cdot \kappa \approx\sum \log q_j = O(s)$ and $\kappa = O(1)$. 
  Computing the product tree for the $q_j$'s and reducing one coefficient of the $u_i$'s modulo each $q_j$
  costs $O\bigl(\Mint{\log Q}\cdot\log k\bigr)$ by \cite[Th.\,10.24]{GG13},
  hence in total $O\bigl(rn\cdot \Mint{s} \log s\bigr)$.
  Using \cref{pr:cplxcrtq}, the loop has a total cost of
  $O\bigl(rn \log Q \cdot ( \Mpol{n}{\kappa} + \log e)\bigr)$.
  By \cite[Th.\,10.25]{GG13}, the final CRT has a negligible $O(n\cdot \Mint{s}\log s)$ cost.
\end{proof}

\begin{remark}
  Likewise, using~\cref{rem:crt-primes}, the overall complexity drops to 
  $O\bigl(rn\Mint{s}\log s + rs\log e\bigr)$
  when many totally split primes exist.
\end{remark}




\subsection{Bad cases: existence of cyclotomic subfield}


In this section, we give the following result, which characterises number fields
which are ``bad'' fields for $e$.
This result can also be found in \cite[\S{}C,~Pr.\,2.1]{Molin}, and we give the proof for completeness.

\begin{theorem}
  \label{thm:equiv-bad-cases}
  Let \(K\) be a number field. The two following assertions are equivalent :
  \begin{enumerate}[\((i)\)]
  \item For almost all prime \(p \in \N\), \(\forall \mf{p} \mid p\),
    \(p^{\rd{\mf p}{p}} \equiv 1 \bmod e\);
  \item \(\Q(\zeta_e)\) is a subfield of \(K\).
  \end{enumerate}
\end{theorem}

For this we will need a result due to Bauer mentioned in~\cite{neukirch}.
\begin{notation}
  Given \(L/K\) a number field extension we denote by \(P(L/K)\) the set
  \( \lbrace \mf p \text{ unramified prime of } K \mid \exists \mf P, \rd{\mf
    P}{\mf p} = 1 \rbrace.  \)
\end{notation}

\begin{lemma}[Bauer in~\cite{neukirch}]
  \label{lem:neuk}
  If \(L / K\) is Galois and \(M / K\) is an arbitrary finite extension, then
  \(
    P(M/K) \dot \subseteq P(L/K)  \iff L \subseteq M.
  \)
\end{lemma}

\begin{proof}[Proof of~\Cref{thm:equiv-bad-cases}]
  The first assertion is true for \(\Q(\zeta_e)\), so \((ii) \implies (i)\) is clear by multiplicativity of the
  inertia degree. Now assume~\((i)\).
  Let \(p \in P(K/\Q)\) and \(\mf p \) s.t.~\(\rd{\mf p}{p} = 1\).
  Condition \((i)\) implies that almost all such prime satisfies
  \(p \equiv 1 \bmod e\), which implies \(p\) is completely split in
  \(\Q(\zeta_e)\), \ie\(p \in P(\Q(\zeta_e)/ \Q)\). Consequently, we have
  \(P( K /\Q) \dot \subseteq P(\Q(\zeta_e) /\Q)\) and \Cref{lem:neuk} gives
  \((ii)\).
\end{proof}

In all generality, the class of fields for which \Cref{alg:crt-root} cannot be
used is larger than the one described by~\Cref{thm:equiv-bad-cases}. However,
if one considers Galois fields only, the two are equivalent.

\subsection{Experimental results}
\label{sec:thome-experimental-results}

\begin{figure*}[htp]
  \subfigure[\(e = 3\)\label{subfig:crt-3}]{\includegraphics[scale=0.193]{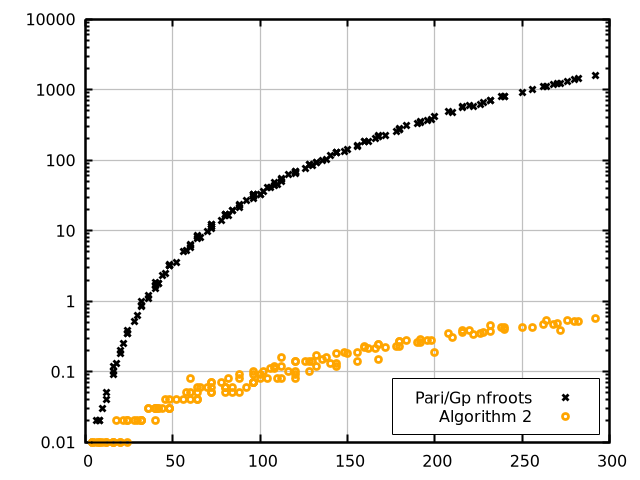}}
  \subfigure[\(e = 71\)\label{subfig:crt-71}]{\includegraphics[scale=0.193]{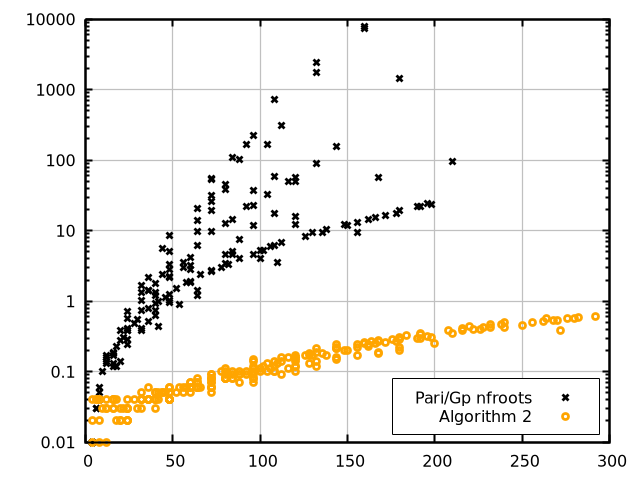}}
  \subfigure[\(e = 1637\)\label{subfig:crt-1637}]{\includegraphics[scale=0.193]{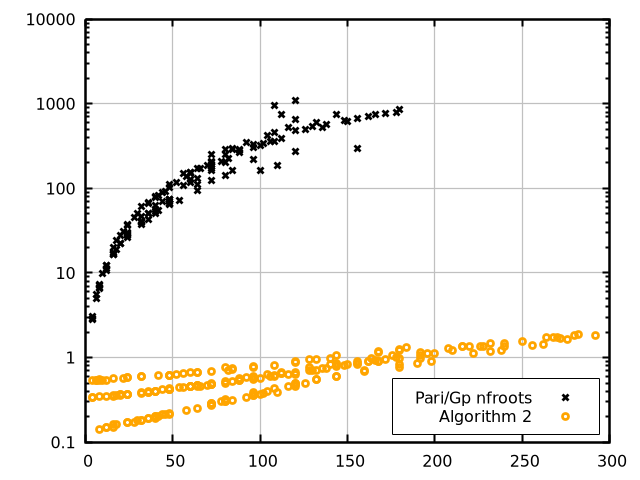}}
  \subfigure[\(e = 13099\)\label{subfig:crt-13099}]{\includegraphics[scale=0.193]{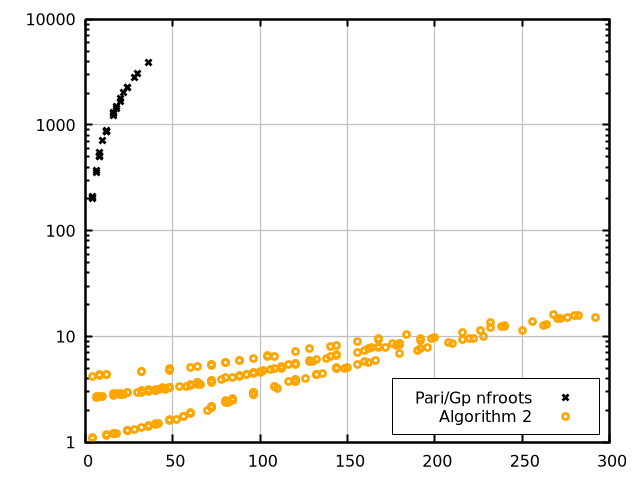}}
  \vspace{-\floatsep}
  \caption{Timings (s) for \texttt{nfroots} and \Cref{alg:crt-root} plotted
    against the degree, for various prime $e$ over cyclotomic fields.}
  \label{fig:thome-naive-root}
\end{figure*}

\begin{table*}[!ht]
  \begin{tabular}{lrrrrrrrrrrrrrr}
    \toprule
    Conductor $m$                           & 113     & 256     & 137     & 149     & 169     & 243     & 249
                                            & 167     & 173     & 179     & 181     & 235     & 197     & 199    \\
    Degree $\varphi(m)$                     & 112     & 128     & 136     & 148     & 156     & 162     & 164
                                            & 166     & 172     & 178     & 180     & 184     & 196     & 198    \\
    $\log_2(e)$, with $e \mid h_{m}^{-}$    & 43.4    & 44.7    & 45.4    & 68.8    & 47.5    & 39.7    & 53.9
                                            & 82.1    & 66.0    & 93.5    & 62.3    & 57.1    & 92.9    & 41.5   \\
    \midrule
    \textit{$\log_2(\Ninf{x})$, $x=y^{1/e}$}& 215.0   & 263.9   & 311.9   & 346.9   & 381.9   & 395.6   & 393.4
                                            & 469.1   & 320.7   & 410.0   & 461.5   & 488.4   & 546.3   & 492.7  \\
    Double CRT (Alg.\,\ref{alg:crt-root})   & 2.5     & 3.4     & 5.7     & 9.3     & 8.7     & 8.6     & 24.7
                                            & 12.3    & 9.1     & 14.0    & 14.6    & 24.5    & 17.3    & 16.2   \\
    \textit{Schirokauer maps}               & \it 176 & \it 15  & \it 748 &\it 1035 &\it 100  &\it 20   &\it 822
                                            &\it 1526 &\it 2408 &\it 4913 &\it 286  &\it 998  &\it 5426 &\it 1864\\
    \bottomrule
  \end{tabular}
  \caption{Timings (s) for \(e\)-th roots within the saturation of Stickelberger \(S\)-units \cite{tw-sti} for
    selected cyclotomic fields}
  \label{tab:timings-s-units}
    
  \vspace{-\floatsep}
\end{table*}

We compared in practice~\cref{alg:crt-root} to standard algorithms
and implementations such as \textsf{Pari/Gp} \texttt{nfroots}. 
For specific exponents
\(e \in \{ 3, 71, 1637, 13099 \}\), and focusing on suitable cyclotomic fields,
we computed the average time taken to compute
the \(e\)-th roots of \(y = x^e\) where \(x\) is a random element with
coefficients of bit size $\log B \in\{1,50,100\}$.
We stress that with this protocol, \cref{alg:crt-root} does \emph{not} take advantage
of being designed for treating factored forms.
The results obtained for $\log B = 100$ are displayed 
in~\cref{fig:thome-naive-root}. 

\begin{remark}
  For small exponents, \ie such that \(3e \leqslant [K:\Q]\), the function
  \textsf{nfroots} from \textsc{Pari/Gp} uses Trager's method~\cite{trager}; 
  in these experiments, this is the case only when $e=3$. Otherwise
  it follows the ideas developed in~\cite{fieker-reconstruction,belabas_roots}.
\end{remark}

From~\cref{fig:thome-naive-root}, we see that our implementation of
\cref{alg:crt-root} using \textsc{SageMath} is, in all cases, much more efficient than
\textsf{nfroots} when the dimension increases. Further, the gap increases with the
exponent. Remark, \eg in \cref{subfig:crt-71}, that \Cref{alg:crt-root} is more
stable than \textsc{Pari/Gp} \textsf{nfroots}.

Finally, as expected, the performances of our algorithm are not much influenced
by the size of the exponent \(e\), to the point where it is perfectly fine to compute $e$-th roots for very large $e$'s. For example, \cref{tab:timings-s-units} shows timings for large $e$'s in a real-life application, namely the computation of $S$-units in cyclotomic fields, by saturating a full-rank multiplicative $S$-units subset arising from Stickelberger's theorem,  as suggested in \cite{tw-sti}.
Concretely, our biggest example here is for $e = 14458667392334948286764635121$ in $\Q(\zeta_{179})$, a 94-bits prime dividing $h_{179}^{-}$, for which an $e$-th root in dimension $n=178$ is computed in only 14 seconds.


\section{A relative Couveignes' method: the bad cases}
\label{sec:roots-bad}

\def\nfsup{K}
\def\nfsub{L}

In this section we will describe how one can compute \(e\)-th roots in the bad
cases, \ie when for all primes \(p \in \N \), there exist at least one
$\mathfrak{p} \mid p$ such that $p^{f_{\mathfrak{p}}} \equiv 1 \mod e$.

\subsection{A relative Couveignes' method}
\label{sec:relat-couv}

Couveignes' method for square-root computation~\cite{couveignes,thome} allows
identifying together the roots modulo a set of inert primes
\(p_1, \dots, p_k\), assuming that the degree of the number field \([K:\Q]\) is
odd. In the following, we show how to generalise this method to
\(e \geqslant 3 \).

A key ingredient of Couveignes' method is the fact that the norm maps
\( \norm{K}{} : K \to \Q\) and
\(\norm{\F_{p^n}}{} : \faktor{\O{K}}{(p)} \simeq \F_{p^n} \to \F_p \) are
coherent when \(p\) is an inert prime integer, meaning that for any \(x \in K\)
we have \(\norm{K} (x) \bmod p = \norm{\F_{p^n}}{} ( x \bmod p) \). We will use
a generalisation of this property to relative extensions.

\begin{lemma}
  \label{lem:norm-com-inert}
  Consider \(\nfsup/\nfsub\) a number field extension, \(\mf{p}\) an inert prime ideal of
  \(\O{\nfsub}\) and \(\mf{P}\) the prime ideal of \(\O{\nfsup}\) above \(\mf{p}\). Then
  the following diagram is commutative.
  \begin{center}
    \begin{tikzcd}[arrow style = tikz]
      \O{\nfsup} \arrow[twoheadrightarrow, r, ""{name=A}] \arrow[d, swap, "\norm{\nfsup/\nfsub}"] 
      & \faktor{\O{\nfsup}}{\mf{P}} \arrow[d, "\norm{}"] \\
      \O{\nfsub} \arrow[twoheadrightarrow, r, ""{name=B}] 
      &  \faktor{\O{\nfsub}}{ \mf{p}}
      \arrow[from=B, to=A, pos=.45, "\scalebox{2.2}{\(\circlearrowleft\)}", phantom]
    \end{tikzcd}
  \end{center}
\end{lemma}

\begin{proof}

  Elements of \(\O {\nfsup}\) can be seen as polynomials with coefficients in \(\nfsub\) with degree
  less than \([\nfsup:\nfsub]\), as 
  \(\O \nfsup \cong \O \nfsub[T] / (P_{\nfsup/\nfsub}(T))\) for some irreducible polynomial
  \(P_{\nfsup/\nfsub} \in \nfsub[T]\).
  Similarly, elements of \(\F_{\mf P} := \smash{\faktor{ \O{\nfsup}}{ \mf P}}\) are
  class of polynomials in \(\F_{\mf p} := \smash[t]{\faktor{\O \nfsub}{\mf p}}\).
  Since \(\mf P \mid \mf p\) is inert, one can define \(\F_{\mf P} / \F_{\mf p}\)
  using the same polynomial \(P_{\nfsup/\nfsub}\), \ie
  \(\F_{\mf P} \cong \F_{\mf p}[T] / (P_{\nfsup/\nfsub}(T))\). 
  Consequently, if \(\alpha\in \O \nfsup\)
  is written as \(\smash[b]{\sum_{i=0}^{[\nfsup:\nfsub]-1}\alpha_i T^i}\), with
  \(\alpha_i \in \O \nfsub\), 
  then its embedding
  into \(\F_{\mf P}\) is expressed as \({\sum_{i=0}^{[\nfsup:\nfsub]-1} \overline{\alpha}_i T^i}\) where
  \(\overline {\cdot} : \O \nfsub \to \F_{\mf p}\).
  Now, recall that the relative norm
  of an element \(\alpha\) in a field extension \(E/F\) is the 
  determinant of the multiplication map \([\alpha]\), seen as a \(F\)-linear
  map of~\(E\).  Since both extensions \(\nfsup/\nfsub\) and \(\F_{\mf P} / \F_{\mf p}\)
  are given by the same defining polynomial, one can consider the action
  \([\alpha]\) and \([\overline \alpha]\) over the ``same'' basis, \ie 
  \([\alpha]\) and \(\overline{\cdot}\) commute. This transfers to the
  determinant, so that
  \(\det [\alpha] \bmod \mf p = \det [ \alpha \bmod \mf P]\).
\end{proof}

Another key ingredient of Couveignes' method for $e=2$ is that, as soon as the degree of the field is odd,
\(\norm{}(\pm x) = \pm \norm{}(x)\).
We present an extension of this property to \(e\)-th roots of unity and
relative extensions.

\begin{lemma}
  \label{lem:norm-bij}
  Let \(E/F\) be a field extension of finite degree and \(e \in \N^*\) such that
  \(\zeta_e\)  \(\in F\). Assume additionnally that
  \(\gcd ([E:F], e) = 1\). Then, for any \(y \in (E^*)^e\), the norm map
  \(\norm{E/F}\) induces a bijection between the zeroes \(Z_E(X^e - y)\) in $E$ and
  \(Z_F\bigl(X^e - \norm{E/F}(y)\bigr)\) in $F$.
\end{lemma}

\begin{proof}
  The set \(Z_E(X^e - y)\) is of the form
  \(\lbrace \zeta_e^i x \mid i \in \segment{0}{e-1} \rbrace\) where \(x\) is a
  fixed \(e\)-th root of \(y\) in \(E\). Similarly, we have
  \[Z_F\bigl(X^e - \norm{E/F}(y)\bigr) = \lbrace \zeta_e^i \cdot \norm{E/F}(x) \mid i \in
  \segment{0}{e-1} \rbrace. \] 
  Now, since \(\zeta_e \in F\),
  \(
    \norm{E/F}\bigl(\zeta_e^i \cdot x\bigr) = \zeta_e^{i[E:F]}\cdot \norm{E/F}(x)\). 
  Thus, it suffices to prove the result for the \(e\)-th roots of
  unity, which follows from the fact that \(\gcd ([E:F], e) = 1\).
\end{proof}

\begin{remark}
  Note that~\Cref{lem:norm-bij} can be applied indifferently to number field extensions or
  finite field extensions. Since both 
  \(\O \nfsub\) and \(\O \nfsub / \mf p\)
  contain
  a primitive \(e\)-th root of unity,~\Cref{lem:norm-bij} applies to both sides
  of the commutative diagram from~\Cref{lem:norm-com-inert}.
\end{remark}

\begin{algorithm}[htb]
  \caption{Relative Couveignes' method for $e$-th root modulo \(p\)}
  \label{alg:crt-couveignes-one-prime}
  \begin{algorithmic}[1]
    \Require A number field extension $\nfsup/\nfsub$ of degree \mbox{prime to~$e\geq 3$},
    $y\in (\nfsup^*)^e$ in factored form, 
    a fixed $e$-th root \(a = \norm{\nfsup/\nfsub}(y)^{1/e}\), 
    a prime integer \(p\) s.t.~each \(\mf p\) above \(p\) in
     \(\O \nfsub\) is inert in \(\nfsup/\nfsub\).
    \Ensure \(x \equiv y^{1/e} \bmod{(p)}\) in \(\nfsup\).
    \State Compute \(\{ \mf p_1, , \dots, \mf p_{{g}} \mid p \O \nfsub = \prod_i \mf p_i  \} \)
    \State Compute \(\{ \mf P_1, , \dots, \mf P_{{g}} \mid \mf p_i \O \nfsup = \mf P_i \}\)
    \For{\(1 \leq i \leq g\)}
    \State \(x_i \gets y^{1/e} \mod \mf P_i\)  \Comment{Pick one root in the residue field}
    \State \(a_i \gets a \bmod \mf p_i\)       \Comment{Expected relative norm mod $\mf p_i$}
    \State \(z_i \gets a_i / \norm{\nfsup/\nfsub}(x_i)\) \Comment{$e$-th root of unity mod $\mf p_i$}
    \State \(x_i \gets x_i\cdot z_i^{[\nfsup:\nfsub]^{-1}\bmod e}\)
    \EndFor
    \State  \Return   \( \text{CRT}_L\bigl(\{x_i\}_{i \in \segment{1}{{g}}}, \{\mf P_i\}_{i \in
      \segment{1}{{g}}}\bigr)\).
  \end{algorithmic}
\end{algorithm}

\begin{lemma}
  \label{lem:gener-couv-one-prime}
  Consider \(\nfsup/\nfsub\) a number field extension and \(e \in \N^*\) such that
  \(\zeta_e \in  \nfsub\). Assume additionnally that
  \(\gcd ([\nfsup:\nfsub], e) = 1\).  Then \Cref{alg:crt-couveignes-one-prime} is correct
  and runs in polynomial time.
\end{lemma}

\begin{proof}
  The output is correct if for each
  \(i \in \segment{1}{g}\), \(x_i\) is the embedding in
  \(\F_{\mf P_i} := \smash[t]{ \faktor{\O{\nfsup}}{\mf P_i}}\) of a fixed \(e\)-th root of \(y\) in \(\nfsup\),
  denoted by \(x\). This is ensured by the combination of
  \cref{lem:norm-com-inert,lem:norm-bij}, which state that in each
  residue field \(\F_{\mf P_i}\) there is exactly one element of 
  \(Z_{\F_{\mf P_i}}(X^e - x_i)\)
  whose relative norm over \(\smash[t]{\faktor{\O \nfsub}{\mf p_i}}\) is \(a_i\).
\end{proof}

Fixing a root in $\nfsub$ of the relative norm of $y$ allows us to make
a consistent choice of residues in $L$ modulo each $\mf P$ from their relative norms
modulo $\mf p$ in $\nfsub$. As shown by \cref{lem:gener-couv-one-prime}, this
allows us to compute a root \(x \bmod p \in \nfsup\) for a given prime
integer \emph{without} the enumeration process
that would otherwise be needed to guess the correct combination of residues in the CRT.
The same idea applies across several prime integers, see~\cref{alg:crt-couveignes}.

\begin{algorithm}[htb]
  \caption{Relative Couveignes' method for $e$-th root}
  \label{alg:crt-couveignes}
  \begin{algorithmic}[1]
    \Require A number field extension $\nfsup/\nfsub$ of degree \mbox{prime to~$e\geq 3$},
    $y\in (\nfsup^*)^e$ in factored form. 
    \Ensure  \(x \in \nfsup^*\) such that $y = x^e$.
    \State Compute $B$ s.t.~$\Ninf{C(x)}\leq B$ \Comment{Using \cref{lem:boundcoef,lem:boundethroot}}
    \State \(a \gets \norm{\nfsup/\nfsub}(y)^{1/e}\) \Comment{Using \S\ref{sec:generic} or
      this algorithm in $\nfsub$}
    \State Choose prime integers $p_1, \dotsc, p_k$ s.t. each prime
    ideal above \(p_j\) in \(\O \nfsub\) is inert in \(\nfsup/\nfsub\) and~$\prod_j p_j \geq 2B$.
    \For{\(1 \leq j \leq k\)}
    \State \(x_j \gets \mathsf{RelativeCouvMod}_p(e, y, a, p_j) \)
    \Comment{~\Cref{alg:crt-couveignes-one-prime}}
    \EndFor
    \State $x \gets  \text{CRT}_{\Z}(\{x_j\}_{j}, \{p_j\}_{j}$ \Comment{Coefficient by coefficient}
    \State  \Return  $x$ with coefficients mapped in $[-B,B]$.
  \end{algorithmic}
\end{algorithm}


\begin{theorem}
  Consider \(\nfsup/\nfsub\) a number field extension and \(e \geq 3\) such that
    \(\zeta_e \in \nfsub\). Assume also that
    \(\gcd ([\nfsup:\nfsub], e) = 1\).  Then, \cref{alg:crt-couveignes} is correct and
    runs in polynomial time in $[\nfsup:\Q]$ and $e$.
\end{theorem}

Note that for~\Cref{alg:crt-couveignes} to be usable as a generic method over a
fixed field \(\nfsup\) together with an exponent \(e\), it is necessary that there is
a subfield \(\nfsub \hookrightarrow \nfsup\) ensuring the existence of infinitely many
prime integers with the right splitting condition mentioned. If \(\nfsup/\nfsub\) is
Galois, one can deduce from Chebotarev density theorem that this is equivalent
to \(\nfsup/\nfsub\) being cyclic~\cite{cohen-advanced,neukirch}.


One of the most costly steps is the computation of
\(\norm{\nfsup/\nfsub}(y)\), whose cost is mitigated by using a factored
form for $y$. The overall complexity depends on the algorithm
used to compute its \(e\)-th root.

\subsection{A recursive relative Couveignes' method}
\label{sec:recurs-relat-couv}
\Cref{alg:crt-couveignes} can be turned into a recursive algorithm, noting that
one important step is to compute the \(e\)-th root of \(\norm{\nfsup/\nfsub}(y)\).
We will focus on the Galois case. Thus, let us fix a Galois number
field \(\nfsup\) and denote by \(G\) its Galois group.

Recall that we assume that \(\Q(\zeta_{e}) \hookrightarrow \nfsup\),
as we are in the ``bad'' case. Let \(\nfsub\) be a subfield of \(\nfsup\) such that
\Cref{alg:crt-couveignes} can be applied, \ie \(\gcd([\nfsup:\nfsub], e) = 1\),
\(\Gal{\nfsup}{\nfsub}\) is cyclic and $\Q(\zeta_e)\hookrightarrow \nfsub$.
In order to apply~\Cref{alg:crt-couveignes}
recursively for the computation of \(\norm{\nfsup/\nfsub}(y)^{1/e}\), one needs the
existence of a subfield of $\nfsub$ satisfying~the~same
conditions. We thus get~\Cref{prop:recursive-couveignes-applicability}, which
describes Galois extensions for which a recursive version
of~\Cref{alg:crt-couveignes} is applicable.

\begin{proposition}
  \label{prop:recursive-couveignes-applicability}
  Consider \(\nfsup/\nfsub\) a Galois extension of number fields. Then one can apply
  a recursive version of~\Cref{alg:crt-couveignes} with respect to \(\nfsup/\nfsub\) and
  \(e \in \N\) if, and only if, \(\nfsup/\nfsub\) is Abelian and \(\gcd([\nfsup:\nfsub], e) = 1\).
\end{proposition}

\begin{proof}
  If \(\nfsup/\nfsub\) is suitable for a recursive version of~\Cref{alg:crt-couveignes},
  then there is a tower of subfields
  \[
    \nfsub = \nfsub_0  \hookrightarrow \nfsub_{1} \hookrightarrow \dots \hookrightarrow \nfsub_{r-1}
    \hookrightarrow \nfsub_r = \nfsup
  \]
  such that for all \(i \in \segment{0}{r-1}\) the extension \(\nfsub_{i+1}/\nfsub_{i}\) is
  cyclic and \(\gcd ([\nfsub_{i+1}:\nfsub_i], e) = 1\). Since
  \(\Gal{\nfsup}{\nfsub} \cong \bigoplus_{i=0}^{r-1} \Gal{\nfsub_{i+1}}{\nfsub_i} \) and
  \([\nfsup:\nfsub] = \prod_{i=1}^{r-1} [\nfsub_{i+1}:\nfsub_i]\), \(\nfsup/\nfsub\) is Abelian with
  \(\gcd([\nfsup:\nfsub], e) =1\). Conversely, denoting by \(G\) the Galois group of
  \(\nfsup/\nfsub\), if this extension is Abelian then \(G\) admits a cyclic
  refinement~\cite{lang2012algebra}. The part on the dimension is clear as well.
\end{proof}

\subsection{Experimental results}
\label{sec:couveignes-simple-exp}


\begin{figure*}[htb]
  \centering
  \subfigure[\(p = 5\)]{\includegraphics[scale=0.2]{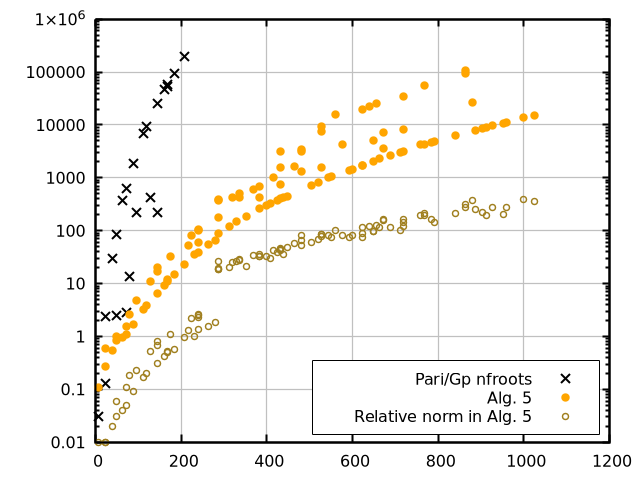}}
  \subfigure[\(p = 11\)]{\includegraphics[scale=0.2]{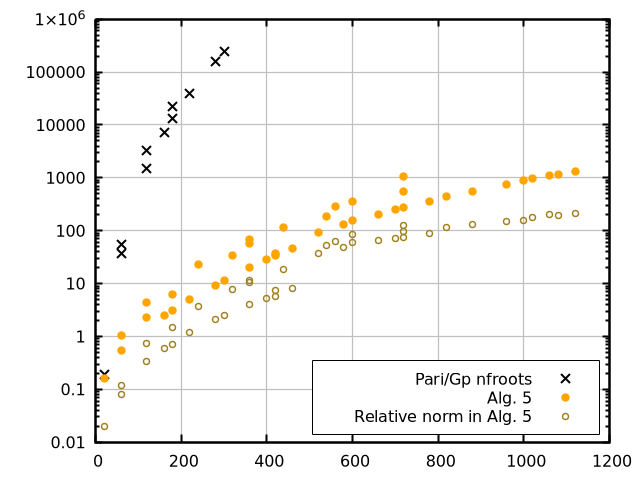}}
  \subfigure[\(p = 23\)]{\includegraphics[scale=0.2]{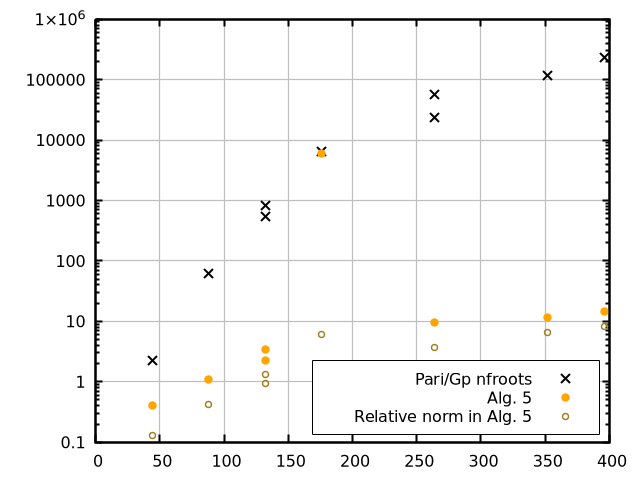}}
    \vspace{-\floatsep}
    \caption{Timings (s) for \texttt{nfroots} and~Alg.~\ref{alg:crt-couveignes}
    plotted against \(n\), over fields \(\Q(\zeta_{pq})\)
    with constant \([\nfsup:\nfsub]= p-1\) and \(e=q\).}
  \label{fig:couveignes-naive-reldeg}
\end{figure*}


\begin{figure*}[htp]
  \centering
  \subfigure[\(p =  5\)]{\includegraphics[scale=0.2]{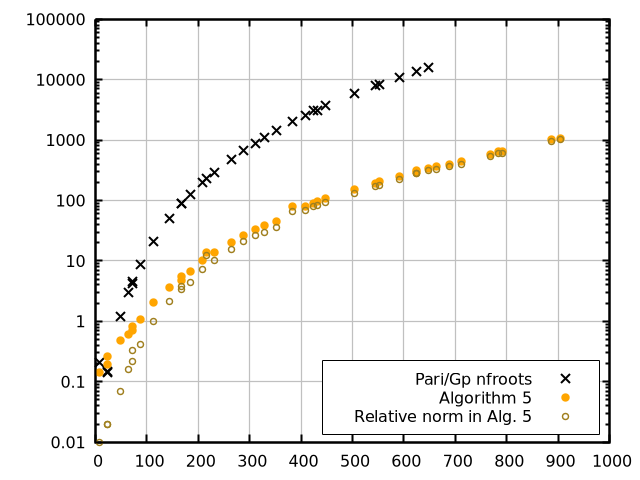}}
  \subfigure[\(p =  11\)]{\includegraphics[scale=0.2]{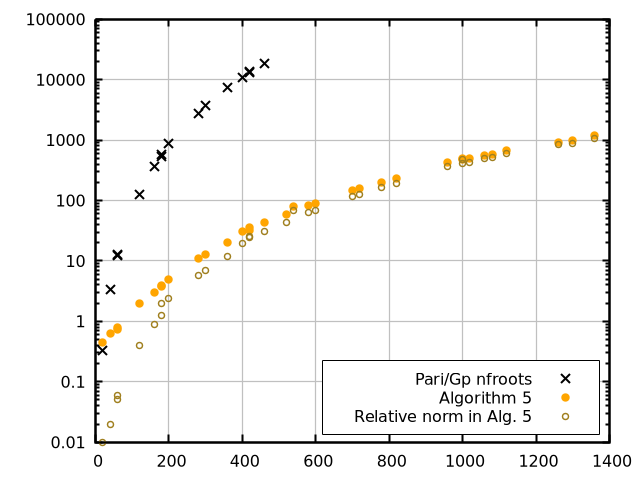}}
  \subfigure[\(p =  23\)]{\includegraphics[scale=0.2]{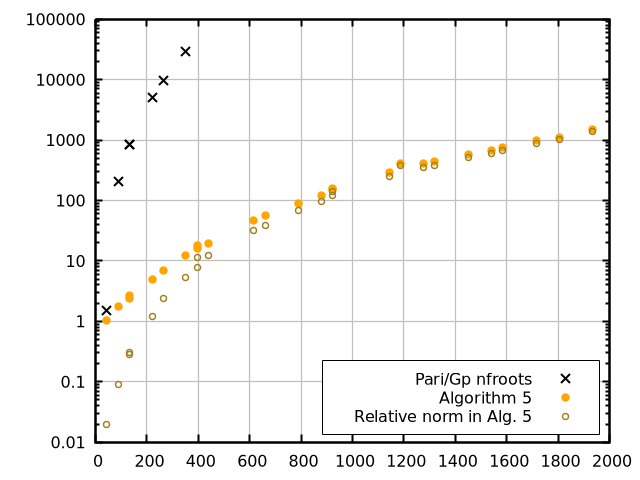}}
  \vspace{-\floatsep}
  \caption{Timings (s) for \texttt{nfroots} and~Alg.~\ref{alg:crt-couveignes}
    plotted against \(n\), over fields \(\Q(\zeta_{pq})\)
    with constant \(e=p\) and \([\nfsup:\nfsub]=\varphi(q)\).}
  \label{fig:couveignes-naive-root}
\end{figure*}


We report now on experimental results we obtained from our implementation of
\Cref{alg:crt-couveignes} in \textsc{SageMath}~\cite{sagemath}. In these
simple experiments we chose to consider cyclotomic fields of the form
\(\Q(\zeta_{pq})\) where \(p\) is a prime integer and \(q \in \N\) satisfies
suitable conditions. The prime \(p\) is constant in each experiment.

We considered two sets of experiments. In the first, we chose \(q\) to be a prime
integer as well, and fixed \([\nfsup:\nfsub] = p\) thus \(e=q\) and \(\nfsub=\Q(\zeta_q)\). In
the second, \(e = p\) so one consider \(\nfsub = \Q(\zeta_e)\) and \([\nfsup:\nfsub] = q \). We
chose to compare the performances of our implementation to the ones of
\textsf{Pari/Gp} \texttt{nfroots}.
Results can be found
in~\cref{fig:couveignes-naive-reldeg,fig:couveignes-naive-root}.


One can see that~\Cref{alg:crt-couveignes} is more efficient than
\textsf{nfroots} in all cases. Recall that \textsf{Pari/Gp} \texttt{nfroots}
uses Trager's method~\cite{trager} when \(3 e \leqslant [\nfsup:\Q]\), which is
the case in our experiments for which \(e=p\). Otherwise, \textsf{Pari/Gp}
\texttt{nfroots} uses the ideas developped
in~\cite{fieker-reconstruction,belabas_roots}. All of these observations tend to
show that our generalisation of Couveignes' method is more effective than
previous algorithms used to compute \(e\)-th roots.

An important constatation is that most of the running time
of~\cref{alg:crt-couveignes} comes either from the computation of
\(\norm{\nfsup/\nfsub}(y)\) when \([\nfsup:\nfsub]\) is large (as predicted by
the theoretical complexity), see~\Cref{fig:couveignes-naive-root}. In practice
it amounts for more than 90\% of the computation time. Improving the efficiency
of this task would render our relative Couveignes' method (and its
implementation) even more impactful. Indeed we compute these norms as products
of polynomials which might not be the most efficient in all cases. Instead, one
could use a half-gcd version of resultants to be quasi-linear.





\section{Saturation: a real-life example}
\label{sec:saturation}

In this section, we first briefly describe the saturation process that can be
used during the computation of $S$-unit groups, when a subgroup of finite index
is already known \cite{multiquadratics,norm_relations,tw-sti}.  Usually, the
saturation is performed for small primes $e$ dividing this index. We show here
the efficiency of our methods to handle much larger values of prime-power $e$'s.

Assume that we have access to \(E = \lbrace y_1, \dots, y_s \rbrace\) a
generating set of \(H\), a subgroup of a multiplicative group \(G\). To fix
ideas, \(G\) is the group \(\U{K, S}\) of \(S\)-units of a number field \(K\)
for some set $S$ of prime ideals, and \(H\) a full-rank
subgroup~\cite{bernard-prepro,tw-sti}. We also assume elements of \(E\)
are given in \emph{factored representation}, \ie their factorisation on a given
multiplicative basis \(U = \{ u_1, \dots, u_r\}\) is known as:
\begin{equation*}
  \forall i \in \segment{1}{s}, \exists (e_{i,j})_{1\leq j\leq r} \in
  \Z^{r} \text{ s.t. } y_i = \prod_{j=1}^r u_{j}^{e_{i, j}}.
\end{equation*}
The overall saturation mechanism can be summarized as follows:
\begin{enumerate}
\item detect elements of $H \cap (K^*)^e$, as we will see, this is actually the
  most costly part;
\item compute the corresponding $e$-th roots, this is exactly the subject of this paper;
\item compute a multiplicative basis of a subgroup of index divided by $e$; this
  part benefits from being realized only once for all possible $e$'s dividing
  the index of $H$ in $G$, in order to contain the size of the elements, and
  won't be discussed further.
\end{enumerate}
A precise theoretical analysis can be found in \cite{norm_relations}.

\subsection{Detecting \texorpdfstring{\(e\)}{e}-th powers}
Let $e$ be a prime-power. In this section we describe how one can
efficiently detect \(e\)-th powers. At a very high level, it goes as follows:
\begin{enumerate}
\item select characters \(\chi_{\mf Q} : (H, \times) \to (\Z/e\Z,  +)\) such that
  \(y \in H^e \implies \chi_{\mf Q}(y) = 0 \) for sufficiently enough primes
  \(\mf Q\);
\item compute the kernel of of \(\chi : y \in H \mapsto (\chi_{\mf Q}(y))_{\mf{Q}} \).
\end{enumerate}
The characters selected in step (1) have to be numerous enough so that the
morphism \(\chi\) contains non trivial \(e\)-th powers.
Theoretically, the validity of this step is controlled by the Grunwald-Wang
theorem, and a practical instantiation of the problematic cases where the
Grunwald-Wang theorem does not directly apply, can be found in
\cite[\S4.2]{BEFHY22}. We shall only keep in mind that problematic cases only
arise when $e$ is a power of $2$.

\subsubsection{Conditions on the primes}
Prime integers \(q\) below suitable \(\mf Q\) verify some conditions. The
reduction map \(\phi_{\mf Q} : \O K \twoheadrightarrow \O K / \mf Q \) needs to
be extendable to \(H\) which is not included in \(\O K\) in general.  Thus,
\(q\) should not divide the numerator or the denominator of any \(y_i\) for
\(i \in \segment{1}{s}\). We will write \(\phi_{\mf Q}\) as well for this extended
map.

Then recall that we wish to detect non trivial \(e\)-th powers, so the residue
field should contain elements which are not \(e\)-th powers. This is equivalent
to \(Q\equiv 1 \mod e\), where \(Q = \norm{K/\Q}(\mf Q)\).

\subsubsection{Definition of the characters}
Once an element \(y\) has been embedded into a residue field \(\O{} / \mf Q\),
one needs to detect whether \(\phi_{\mf Q}(y)\) is an \(e\)-th power or not. One
can note that for any element \(t \in \F_Q^*\), its power
\(\smash{t^{(Q-1)/e}}\) is an \(e\)-th root of unity, and is equal to \(1\) if, and only
if, \(t \in (\F_Q^*)^e\). Thus one puts
\(\chi_{\mf Q} : y \mapsto \log_{\zeta_e}( \phi_{\mf Q}(y)^{(Q-1)/e}) \), with
\(\zeta_e\) is a primitive \(e\)-th root of unity.

Actually, this infinite family of characters can be completed with easier-to-compute
alternative characters, but of which we only have finitely many.
For instance, when considering multiplicative subsets of $S$-units for some set $S$ of prime ideals,
we can include the $\mf p$-adic valuations modulo $e$, for each $\mf p \in S$.
Another very important finite family of characters is given by \emph{Schirokauer maps} \cite{Schiro}, which can be viewed as an approximation of the $e$-adic logarithm. Namely, in our large $e$'s experiments, we compute
\[
\lambda: y\in K \longmapsto \frac{y^{\rho_e} - 1 \bmod e^2}{e},
\]
where $\rho_e$ is the least common multiple of the $\bigl\lvert(\O K / \mf e)^*\bigr\rvert$ for all prime ideals $\mf e$ dividing $e\O K$. This yields $[K:\Q]$ characters ---~one per coefficient~--- modulo $e$ \cite[Pr.\,3.8]{Schiro}.

\subsubsection{Number of characters} In order to detect non trivial powers,
\ie elements in \(G^e \setminus H^e\), one only has to intersect $\ker \chi_{\mathfrak{Q}}$
for sufficiently many $\mathfrak{Q}$. If $s$ is the cardinal of a generating family
\(E\) of $H$, then the rank of $H/(H \cap (K^*)^e)$ is $s' \leqslant s$. If we
consider the $\chi_{\mathfrak{Q}}$ to be uniformly distributed in the dual
then~\cite[Lem.\,8.2]{number_field_sieve} can be modified to show that $s' + r$
characters generate the dual with probability at least $1-e^{-r}$.

\subsection{Practical considerations}

\subsubsection{Computing suitable primes}
Discarding the condition regarding the denominators of elements of \(E\),
suitable prime integers \(q \mid \mf Q\) only need to verify
\(e \mid \norm{K/\Q}(\mf Q) - 1\). As we have to compute a discrete logarithm in $\F_Q$,
it is desirable to restrict to primes of inertia degree $1$.
Hence, instead of drawing random primes of a given bit-length, it is best to test primality on integers $q \equiv 1\mod e$ until sufficiently many primes of inertia degree $1$ are
found.
Note that, contrary to the CRT case, small primes give as much information as big primes for the characters.

In the Galois case, this comes down to find completely split primes \(q\), and Chebotarev density theorem assures us that the density of those primes in the set of all prime integers is \(1 / [K:\Q]\).
Hence, we expect to find $k \geq s'+r$ suitable primes $q\equiv 1\mod e$ of a given bitsize in \(O(k [K:\Q])\) trials.

For cyclotomic fields  \(K_m = \Q(\zeta_m)\) the situation is even better, as
a prime \(q\) completely
splits in $K$ if, and only if, \(q \equiv 1 \bmod m\).
Thus one needs to find $k$ primes that directly verify the congruence $q \equiv 1\mod \lcm(e,m)$, which can be done at a given bitsize in $O(k)$ trials.


\subsubsection{Computing non trivial relations}
Once characters have been selected as a set of prime ideals
\(S = \{ \mf Q_1, \dots, \mf Q_k \}\) above suitable prime integers \(q\), one needs to
compute the values \(\chi_{\mf Q}(y_i)\) for all~\(\mf Q\) and
\(i \in \segment{1}{s}\), then identify the kernel of \(\chi\) yielding
\(H \cap (K^*)^e\).

Recall that each element \(y_i\) of the generating set \(E\) is known through
its decomposition in the multiplicative basis \(U\). Since a character is a morphism,
the image of \(\chi_{\mf Q}\) is determined by the collection of elements
\(\chi_{\mf Q}(u_j) \in \Z/e\Z\).
As the $u_j$ are expected to be small, it is generally more efficient to compute
first all $\chi_{\mf Q}(u_j)$ and reconstruct each \(\chi_{\mf Q}(y_i)\)
as \(\sum_{j=1}^r e_{i, j} \chi_{\mf Q}(u_j)\).

The image \(\chi(H)\) is then generated by the rows of the matrix
\[
M = \Bigl( \chi_{\mf Q}(y_i) \Bigr)_{\overset{i \in \segment{1}{s}}{\mf Q \in S}} \in
  \mathcal{M}_{s, k}\bigl(\Z/e\Z\bigr).
\]
Therefore, the left kernel of \(M\) in \(\Z/e\Z\) yields
(after lifting to $\Z$) vectors
\((\alpha_1, \dots, \alpha_t) \in \segment{0}{e-1}^t\)
such that \(\prod_{i=1}^s y_i^{\alpha_i} \in (K^*)^e \).

Note that when $e$ is a prime-power, we can use Howell's normal form
\cite{Storjohann}. This form might give redundant solutions, but those will be
sorted out in the final reconstruction phase.

\subsubsection{Complexity analysis} Actually, detecting of \(e\)-th powers can
be problematically long if the sizes of the different finite fields are too
large. This mecanically happens when \(e\) is large. The bottleneck is the
computation of discrete logarithms in subgroups of order~\(e\), and we need to
do it \(rk\) times. The cost is then \(O(rk \sqrt{e})\) using generic
techniques. It can be mitigated by using index calculus methods, which can also
benefit from a precomputation phase when computing many ($r$ here) discrete logarithms.

Instead, in our practical case of interest, $\mf p$-adic valuations
and Schirokauer maps give just enough characters to be able to detect $e$-th powers.
Considering that $\mf p$-adic valuations are known (this is arguably often the case),
the cost of detecting $e$-th powers amounts to computing those Schirokauer maps.
In cyclotomic fields, $\rho_e$ is at most $e^n-1$, so the total cost of detection is
\(
O (r \cdot n \log e \cdot \Mpol{n}{2 \log e}).
\)
Timings for the computation of Schirokauer maps are given in \cref{tab:timings-s-units}, 
and unexpected variations directly relate to the size of $\rho_e$.

\subsection{Experiments}
\begin{figure*}[htp]
  \centering \subfigure[\(e \not \mid m\)]
  {\includegraphics[scale=0.2]{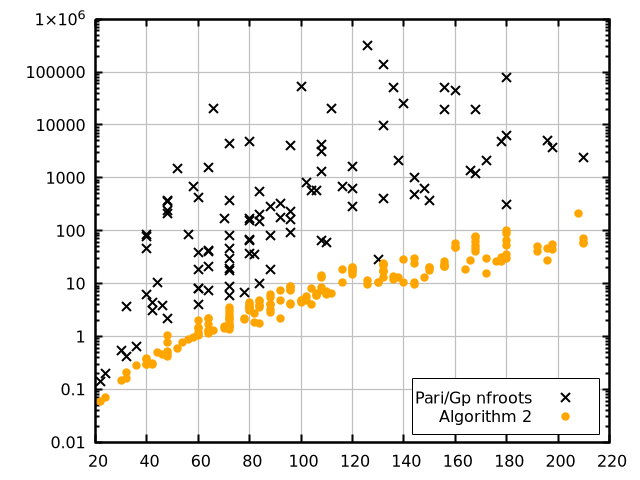}}\qquad
  \subfigure[\(e \mid m\)]{\includegraphics[scale=0.2]{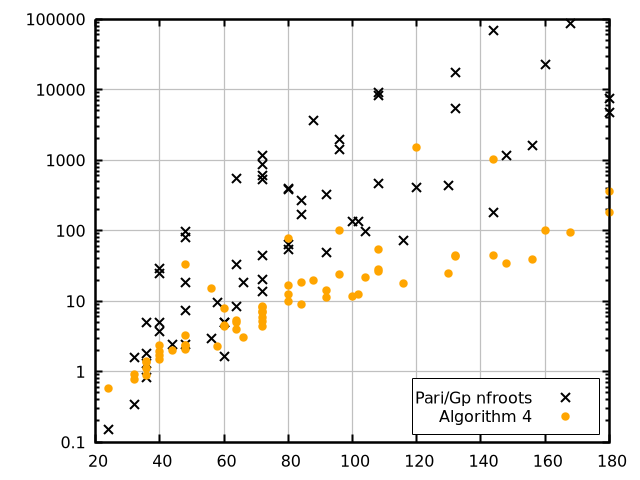}}
  \vspace{-\floatsep}
  \caption{Timings (s) for \texttt{nfroots} and Alg. 4 plotted against
    the dimension for saturation process.}
  \label{fig:couveignes-saturation}
\end{figure*}

In this section, we report the impact of the algorithms presented here on the
running time in saturation processes while computing \(S\)-units in the manner
of~\cite{tw-sti}. In particular, elements are always given in factored
form, 
and we do \emph{never} require to compute the product.
Timings can
be found in~\Cref{fig:couveignes-saturation} and \cref{tab:timings-s-units}.
We separated ``good'' cases where
\(e \not \mid m\) and ``bad'' cases where \(\Q(\zeta_e) \subseteq K\), and we selected
\(e\) to be the largest prime factor of \(h_m^-\) (resp. s.t. \(\gcd(m, h_m^-)\)).

We note a significant gain using~\Cref{alg:crt-root} in good
cases, especially when \(n\) is large. In bad cases, \Cref{alg:crt-couveignes}
also outperforms \textsc{Pari/Gp} \textsf{nfroots} when \(n\) is large,
though its impact is smaller in this case.
The running time of \textsf{nfroots} is greatly influenced
by the size of \(e\), whereas~\Cref{alg:crt-root} is relatively stable with
respect to this parameter. Finally, we stress that much data is missing for
\textsc{Pari/Gp} \textsf{nfroots}, due to its asymptotical limitations.

In order to get a more precise
picture of the power of the Double-CRT algorithm, we gathered
in \Cref{tab:timings-s-units} its performances for very large exponents \(e\).




\bibliographystyle{ACM-Reference-Format}
\bibliography{biblio.bib}


\newcommand{\noopsort}[1]{}
\begin{thebibliography}{25}


\ifx \showCODEN    \undefined \def \showCODEN     #1{\unskip}     \fi
\ifx \showDOI      \undefined \def \showDOI       #1{#1}\fi
\ifx \showISBNx    \undefined \def \showISBNx     #1{\unskip}     \fi
\ifx \showISBNxiii \undefined \def \showISBNxiii  #1{\unskip}     \fi
\ifx \showISSN     \undefined \def \showISSN      #1{\unskip}     \fi
\ifx \showLCCN     \undefined \def \showLCCN      #1{\unskip}     \fi
\ifx \shownote     \undefined \def \shownote      #1{#1}          \fi
\ifx \showarticletitle \undefined \def \showarticletitle #1{#1}   \fi
\ifx \showURL      \undefined \def \showURL       {\relax}        \fi
\providecommand\bibfield[2]{#2}
\providecommand\bibinfo[2]{#2}
\providecommand\natexlab[1]{#1}
\providecommand\showeprint[2][]{arXiv:#2}

\bibitem[\protect\citeauthoryear{Adleman, Manders, and Miller}{Adleman
  et~al\mbox{.}}{1977}]%
        {AMM77}
\bibfield{author}{\bibinfo{person}{Leonard Adleman}, \bibinfo{person}{Kenneth
  Manders}, {and} \bibinfo{person}{Gary~L. Miller}.}
  \bibinfo{year}{1977}\natexlab{}.
\newblock \showarticletitle{On Taking Roots in Finite Fields}. In
  \bibinfo{booktitle}{\emph{18th Annual Symposium on Foundations of Computer
  Science}}. \bibinfo{publisher}{IEEE Computer Society},
  \bibinfo{pages}{175--178}.
\newblock


\bibitem[\protect\citeauthoryear{Bauch, Bernstein, de~Valence, Lange, and van
  Vredendaal}{Bauch et~al\mbox{.}}{2017}]%
        {multiquadratics}
\bibfield{author}{\bibinfo{person}{Jens Bauch}, \bibinfo{person}{Daniel~J.
  Bernstein}, \bibinfo{person}{Henry de Valence}, \bibinfo{person}{Tanja
  Lange}, {and} \bibinfo{person}{Christine van Vredendaal}.}
  \bibinfo{year}{2017}\natexlab{}.
\newblock \showarticletitle{{Short Generators Without Quantum Computers: The
  Case of Multiquadratics}}. In \bibinfo{booktitle}{\emph{Advances in
  Cryptology -- EUROCRYPT 2017}},
  \bibfield{editor}{\bibinfo{person}{Jean-S{\'e}bastien Coron} {and}
  \bibinfo{person}{Jesper~Buus Nielsen}} (Eds.). \bibinfo{publisher}{Springer
  International Publishing}, \bibinfo{address}{Cham}, \bibinfo{pages}{27--59}.
\newblock
\showISBNx{978-3-319-56620-7}


\bibitem[\protect\citeauthoryear{Belabas}{Belabas}{2004}]%
        {belabas_roots}
\bibfield{author}{\bibinfo{person}{Karim Belabas}.}
  \bibinfo{year}{2004}\natexlab{}.
\newblock \showarticletitle{{A relative van Hoeij algorithm over number
  fields}}.
\newblock \bibinfo{journal}{\emph{J. Symb. Comput.}}  \bibinfo{volume}{37}
  (\bibinfo{date}{05} \bibinfo{year}{2004}), \bibinfo{pages}{641--668}.
\newblock
\urldef\tempurl%
\url{https://doi.org/10.1016/j.jsc.2003.09.003}
\showDOI{\tempurl}


\bibitem[\protect\citeauthoryear{Bernard, Lesavourey, Nguyen, and
  Roux-Langlois}{Bernard et~al\mbox{.}}{2022}]%
        {tw-sti}
\bibfield{author}{\bibinfo{person}{Olivier Bernard}, \bibinfo{person}{Andrea
  Lesavourey}, \bibinfo{person}{Tuong-Huy Nguyen}, {and}
  \bibinfo{person}{Adeline Roux-Langlois}.} \bibinfo{year}{2022}\natexlab{}.
\newblock \showarticletitle{Log-$\mathcal{S}$-unit lattices using Explicit
  Stickelberger Generators to solve Approx Ideal-SVP}. In
  \bibinfo{booktitle}{\emph{Advances in Cryptology -- ASIACRYPT 2022}}
  \emph{(\bibinfo{series}{LNCS}, Vol.~\bibinfo{volume}{13793})}.
  \bibinfo{publisher}{Springer}, \bibinfo{pages}{677--708}.
\newblock


\bibitem[\protect\citeauthoryear{Bernard and Roux-Langlois}{Bernard and
  Roux-Langlois}{2020}]%
        {bernard-prepro}
\bibfield{author}{\bibinfo{person}{Olivier Bernard} {and}
  \bibinfo{person}{Adeline Roux-Langlois}.} \bibinfo{year}{2020}\natexlab{}.
\newblock \showarticletitle{Twisted-PHS: Using the Product Formula to Solve
  Approx-SVP in Ideal Lattices}. In \bibinfo{booktitle}{\emph{Advances in
  Cryptology -- ASIACRYPT 2020}}, \bibfield{editor}{\bibinfo{person}{Shiho
  Moriai} {and} \bibinfo{person}{Huaxiong Wang}} (Eds.).
  \bibinfo{publisher}{Springer}, \bibinfo{pages}{349--380}.
\newblock


\bibitem[\protect\citeauthoryear{Biasse, Erukulangara, Fieker, Hofmann, and
  Youmans}{Biasse et~al\mbox{.}}{2022}]%
        {BEFHY22}
\bibfield{author}{\bibinfo{person}{Jean-François Biasse},
  \bibinfo{person}{Muhammed~R. Erukulangara}, \bibinfo{person}{Claus Fieker},
  \bibinfo{person}{Tommy Hofmann}, {and} \bibinfo{person}{William Youmans}.}
  \bibinfo{year}{2022}\natexlab{}.
\newblock \showarticletitle{Mildly Short Vectors in Ideals of Cyclotomic Fields
  Without Quantum Computers}.
\newblock \bibinfo{journal}{\emph{Mathematical Cryptology}}
  \bibinfo{volume}{2}, \bibinfo{number}{1} (\bibinfo{date}{Nov.}
  \bibinfo{year}{2022}), \bibinfo{pages}{84–--107}.
\newblock
\urldef\tempurl%
\url{https://journals.flvc.org/mathcryptology/article/view/132573}
\showURL{%
\tempurl}


\bibitem[\protect\citeauthoryear{Biasse, Fieker, Hofmann, and Page}{Biasse
  et~al\mbox{.}}{2020}]%
        {norm_relations}
\bibfield{author}{\bibinfo{person}{Jean-François Biasse},
  \bibinfo{person}{Claus Fieker}, \bibinfo{person}{Tommy Hofmann}, {and}
  \bibinfo{person}{Aurel Page}.} \bibinfo{year}{2020}\natexlab{}.
\newblock \bibinfo{title}{Norm relations and computational problems in number
  fields}.
\newblock
\newblock
\showeprint[arxiv]{2002.12332}~[math.NT]


\bibitem[\protect\citeauthoryear{Brent}{Brent}{1975}]%
        {Bre75}
\bibfield{author}{\bibinfo{person}{Richard~P. Brent}.}
  \bibinfo{year}{1975}\natexlab{}.
\newblock \showarticletitle{{Multiple-precision zero-finding methods and the
  complexity of elementary function evaluation}}.
\newblock \bibinfo{journal}{\emph{Analytic Computational Complexity}}
  (\bibinfo{year}{1975}), \bibinfo{pages}{151--176}.
\newblock


\bibitem[\protect\citeauthoryear{Buhler, Lenstra, and Pomerance}{Buhler
  et~al\mbox{.}}{1993}]%
        {number_field_sieve}
\bibfield{author}{\bibinfo{person}{Joe~P. Buhler}, \bibinfo{person}{Hendrik~W.
  Lenstra}, {and} \bibinfo{person}{Carl Pomerance}.}
  \bibinfo{year}{1993}\natexlab{}.
\newblock \showarticletitle{Factoring integers with the number field sieve}. In
  \bibinfo{booktitle}{\emph{The development of the number field sieve}},
  \bibfield{editor}{\bibinfo{person}{Arjen~K. Lenstra} {and}
  \bibinfo{person}{Hendrik~W. Lenstra}} (Eds.). \bibinfo{publisher}{Springer},
  \bibinfo{pages}{50--94}.
\newblock
\showISBNx{978-3-540-47892-8}


\bibitem[\protect\citeauthoryear{Cohen}{Cohen}{1993}]%
        {cohen-course}
\bibfield{author}{\bibinfo{person}{Henri Cohen}.}
  \bibinfo{year}{1993}\natexlab{}.
\newblock \bibinfo{booktitle}{\emph{A Course in Computational Algebraic Number
  Theory}}.
\newblock \bibinfo{publisher}{Springer}.
\newblock
\showISBNx{0-387-55640-0}


\bibitem[\protect\citeauthoryear{Cohen}{Cohen}{2012}]%
        {cohen-advanced}
\bibfield{author}{\bibinfo{person}{Henri Cohen}.}
  \bibinfo{year}{2012}\natexlab{}.
\newblock \bibinfo{booktitle}{\emph{Advanced Topics in Computational Number
  Theory}}.
\newblock \bibinfo{publisher}{Springer}.
\newblock
\showISBNx{9781441984890}


\bibitem[\protect\citeauthoryear{Couveignes}{Couveignes}{1997}]%
        {couveignes}
\bibfield{author}{\bibinfo{person}{Jean-Marc Couveignes}.}
  \bibinfo{year}{1997}\natexlab{}.
\newblock \showarticletitle{Computing A Square Root For The Number Field
  Sieve}.
\newblock   \bibinfo{volume}{1554} (\bibinfo{date}{06} \bibinfo{year}{1997}).
\newblock
\showISBNx{978-3-540-57013-4}
\urldef\tempurl%
\url{https://doi.org/10.1007/BFb0091540}
\showDOI{\tempurl}


\bibitem[\protect\citeauthoryear{Fieker and Friedrichs}{Fieker and
  Friedrichs}{2000}]%
        {fieker-reconstruction}
\bibfield{author}{\bibinfo{person}{Claus Fieker} {and} \bibinfo{person}{Carsten
  Friedrichs}.} \bibinfo{year}{2000}\natexlab{}.
\newblock \showarticletitle{On Reconstruction of Algebraic Numbers}. In
  \bibinfo{booktitle}{\emph{Algorithmic Number Theory}},
  \bibfield{editor}{\bibinfo{person}{Wieb Bosma}} (Ed.).
  \bibinfo{publisher}{Springer}, \bibinfo{pages}{285--296}.
\newblock
\showISBNx{978-3-540-44994-2}


\bibitem[\protect\citeauthoryear{Gama and Nguyen}{Gama and Nguyen}{2008}]%
        {GN08}
\bibfield{author}{\bibinfo{person}{Nicolas Gama} {and}
  \bibinfo{person}{Phong~Q. Nguyen}.} \bibinfo{year}{2008}\natexlab{}.
\newblock \showarticletitle{Predicting Lattice Reduction}. In
  \bibinfo{booktitle}{\emph{{EUROCRYPT}}} \emph{(\bibinfo{series}{LNCS},
  Vol.~\bibinfo{volume}{4965})}. \bibinfo{publisher}{Springer},
  \bibinfo{pages}{31--51}.
\newblock


\bibitem[\protect\citeauthoryear{{Joachim von zur Gathen and J{\"{u}}rgen
  Gerhard}}{{Joachim von zur Gathen and J{\"{u}}rgen Gerhard}}{2013}]%
        {GG13}
\bibfield{author}{\bibinfo{person}{{Joachim von zur Gathen and J{\"{u}}rgen
  Gerhard}}.} \bibinfo{year}{2013}\natexlab{}.
\newblock \bibinfo{booktitle}{\emph{Modern Computer Algebra}
  (\bibinfo{edition}{3} ed.)}.
\newblock \bibinfo{publisher}{Cambridge University Press}.
\newblock


\bibitem[\protect\citeauthoryear{Lang}{Lang}{2012}]%
        {lang2012algebra}
\bibfield{author}{\bibinfo{person}{Serge Lang}.}
  \bibinfo{year}{2012}\natexlab{}.
\newblock \bibinfo{booktitle}{\emph{Algebra}}. Vol.~\bibinfo{volume}{211}.
\newblock \bibinfo{publisher}{Springer Science \& Business Media}.
\newblock


\bibitem[\protect\citeauthoryear{Molin}{Molin}{2010}]%
        {Molin}
\bibfield{author}{\bibinfo{person}{Pascal Molin}.}
  \bibinfo{year}{2010}\natexlab{}.
\newblock \emph{\bibinfo{title}{Int\'{e}gration num\'erique et calculs de
  fonctions {L}}}.
\newblock \bibinfo{thesistype}{Ph.D. Dissertation}.
  \bibinfo{school}{L'Universit\'e Bordeaux I}.
\newblock


\bibitem[\protect\citeauthoryear{Neukirch}{Neukirch}{1999}]%
        {neukirch}
\bibfield{author}{\bibinfo{person}{J{\"{u}}rgen Neukirch}.}
  \bibinfo{year}{1999}\natexlab{}.
\newblock \bibinfo{booktitle}{\emph{Algebraic Number Theory}}.
  \bibinfo{series}{Grundlehren der mathematischen Wissenschaften},
  Vol.~\bibinfo{volume}{322}.
\newblock \bibinfo{publisher}{Springer Berlin, Heidelberg}.
\newblock


\bibitem[\protect\citeauthoryear{\noopsort{Pari}{The
  PARI~Group}}{\noopsort{Pari}{The PARI~Group}}{2022}]%
        {PARI2}
\noopsort{Pari}{The PARI~Group} \bibinfo{year}{2022}\natexlab{}.
\newblock \bibinfo{booktitle}{\emph{{PARI/GP version \texttt{2.13.4}}}}.
\newblock \noopsort{Pari}{The PARI~Group}, Univ. Bordeaux.
\newblock
\newblock
\shownote{available from \url{http://pari.math.u-bordeaux.fr/}.}


\bibitem[\protect\citeauthoryear{\noopsort{Sage}{The Sage
  Developers}}{\noopsort{Sage}{The Sage Developers}}{2023}]%
        {sagemath}
\bibfield{author}{\bibinfo{person}{\noopsort{Sage}{The Sage Developers}}.}
  \bibinfo{year}{2023}\natexlab{}.
\newblock \bibinfo{booktitle}{\emph{{S}ageMath, the {S}age {M}athematics
  {S}oftware {S}ystem ({V}ersion x.y.z)}}.
\newblock
\newblock
\shownote{\url{https://www.sagemath.org}.}


\bibitem[\protect\citeauthoryear{Schirokauer}{Schirokauer}{1993}]%
        {Schiro}
\bibfield{author}{\bibinfo{person}{Oliver Schirokauer}.}
  \bibinfo{year}{1993}\natexlab{}.
\newblock \showarticletitle{Discrete logarithms and local units}.
\newblock \bibinfo{journal}{\emph{Philosophical Transactions: Physical Sciences
  and Engineering}} \bibinfo{volume}{345}, \bibinfo{number}{1676}
  (\bibinfo{year}{1993}), \bibinfo{pages}{409--423}.
\newblock


\bibitem[\protect\citeauthoryear{Shoup}{Shoup}{1993}]%
        {Sho93}
\bibfield{author}{\bibinfo{person}{Victor Shoup}.}
  \bibinfo{year}{1993}\natexlab{}.
\newblock \showarticletitle{Factoring polynomials over finite fields:
  Asymptotic complexity vs.~reality}. In \bibinfo{booktitle}{\emph{Proceedings
  of the IMACS Symposium}}. \bibinfo{pages}{124--129}.
\newblock


\bibitem[\protect\citeauthoryear{Storjohann}{Storjohann}{2013}]%
        {Storjohann}
\bibfield{author}{\bibinfo{person}{Arne Storjohann}.}
  \bibinfo{year}{2013}\natexlab{}.
\newblock \emph{\bibinfo{title}{Algorithms for Matrix Canonical Forms}}.
\newblock \bibinfo{thesistype}{Ph.D. Dissertation}. \bibinfo{school}{Swiss
  Federal Institute of Technology, Zurich}.
\newblock


\bibitem[\protect\citeauthoryear{Thom{\'e}}{Thom{\'e}}{2012}]%
        {thome}
\bibfield{author}{\bibinfo{person}{Emmanuel Thom{\'e}}.}
  \bibinfo{year}{2012}\natexlab{}.
\newblock \showarticletitle{Square Root Algorithms for the Number Field Sieve}.
  In \bibinfo{booktitle}{\emph{Arithmetic of Finite Fields}},
  \bibfield{editor}{\bibinfo{person}{Ferruh {\"O}zbudak} {and}
  \bibinfo{person}{Francisco Rodr{\'i}guez-Henr{\'i}quez}} (Eds.).
  \bibinfo{publisher}{Springer}, \bibinfo{pages}{208--224}.
\newblock
\showISBNx{978-3-642-31662-3}


\bibitem[\protect\citeauthoryear{Trager}{Trager}{1976}]%
        {trager}
\bibfield{author}{\bibinfo{person}{Barry~M. Trager}.}
  \bibinfo{year}{1976}\natexlab{}.
\newblock \showarticletitle{Algebraic Factoring and Rational Function
  Integration}. In \bibinfo{booktitle}{\emph{Proceedings of the Third ACM
  Symposium on Symbolic and Algebraic Computation}}
  \emph{(\bibinfo{series}{SYMSAC '76})}. \bibinfo{publisher}{ACM},
  \bibinfo{pages}{219–226}.
\newblock
\showISBNx{9781450377904}


\end{thebibliography}


\end{document}